\documentclass[final,1p,times]{elsarticle}
\usepackage{lineno,hyperref}
\usepackage{amsmath,amsfonts}
\usepackage{verbatim}
\usepackage{latexsym}
\usepackage{graphicx}
\usepackage{epstopdf}
\usepackage{subfig}
\usepackage{caption}
\usepackage{dsfont}
\usepackage{indentfirst}
\usepackage{xcolor}
\usepackage[notref,notcite]{showkeys}
\usepackage{algorithm}
\usepackage{algpseudocode}
\usepackage{booktabs}
\usepackage{natbib}

\allowdisplaybreaks[3]
\newtheorem{theorem}{Theorem}[section]
\newtheorem{lemma}[theorem]{Lemma}
\newtheorem{corollary}[theorem]{Corollary}

\newtheorem{assumption}[theorem]{Assumption}
\newtheorem{example}[theorem]{Example}

\newtheorem{remark}[theorem]{Remark}

\newcommand{\tr}{\tau_n }
\newcommand{\eds}{e_{\Delta}(s^-) }
\newcommand{\tdr}{\theta_{\Delta,n} }
\newcommand{\rdr}{\rho_{\Delta,n} }
\newcommand{\muh}{\mu ^ {-1} (\varphi(\Delta ))}
\newcommand{\xx}{|x_0|^2}
\newcommand{\intt}{\int_0^t}
\newcommand{\ve}{\vee}
\newcommand{\we}{\wedge}
\newcommand{\Dst}{\Delta ^ {\ast}}
\newcommand{\qu}{\quad}
\newcommand{\no}{\nonumber \\}
\newcommand{\bxdt}{\bar x_{\Delta }(t^-)}
\newcommand{\bxds}{\bar x_{\Delta }(s^-)}
\newcommand{\xdt}{ x_{\Delta }(t^-)}
\newcommand{\xds}{ x_{\Delta }(s^-)}
\newcommand{\XD}{X_{\Delta}}
\newcommand{\xdtk}{X_{\Delta} (t_k^-)}

\def\a{\alpha} \def\g{\gamma}  
\def\b{\beta} \def\veps{\varepsilon}
\def\l{\left}  \def\r{\right}   \def\B{\Big} \def\ph{\varphi}

\newcommand{\tkk}{t_{k+1}^-}
\newcommand{\tk}{t_{k}^-}
\newcommand{\fdxs}{f_{\Delta} (\bar x _{\Delta}(s^-)) }
\newcommand{\gdxs}{g_{\Delta} (\bar x _{\Delta}(s^-)) }
\newcommand{\hdxs}{h_{\Delta} (\bar x _{\Delta}(s^-)) }

\newcommand{\fdxt}{f_{\Delta} (\bar x _{\Delta}(t^-)) }
\newcommand{\gdxt}{g_{\Delta} (\bar x _{\Delta}(t^-)) }

\newcommand{\fD}{f_{\Delta}}
\newcommand{\gD}{g_{\Delta}}
\newcommand{\hD}{h_{\Delta}}

\newcommand{\DD}{\Delta}
\newcommand{\la}{\lambda}
\newcommand{\Ito}{It\^{o} formula }
\newcommand{\Holder}{H\"{o}lder inequality }

\newcommand{\As}{Assumption }
\newcommand{\pr}{\textbf{Proof.} }

\newcommand{\E}{\mathbb{E}}
\newcommand{\PP}{\mathbb{P}}
\newcommand{\II}{\mathbb{I}}

\newcommand{\R}{\mathbb{R}}

\makeatletter
\@addtoreset{equation}{section}
\makeatother

\modulolinenumbers[5]

\begin{document}

\begin{frontmatter}
\title{The truncated EM method for stochastic differential equations with Poisson jumps}


\author[mainaddress]{Shounian Deng}

\author[secondaryaddress]{Weiyin Fei \corref{correspondingauthor}}
\cortext[correspondingauthor]{Corresponding author}
\ead{wyfei@ahpu.edu.cn}

\author[address3]{Wei Liu }
\author[thirdaryaddress]{Xuerong Mao }
\address[mainaddress]{School of Science, Nanjing University of Science and Technology, Nanjing, Jiangsu 210094, China}
\address[secondaryaddress]{School of Mathematics and Physics, Anhui Polytechnic University, Wuhu, Anhui 24100, China}
\address[address3]{Department of Mathematics, Shanghai Normal University, Shanghai 200234, China}
\address[thirdaryaddress]{Department of Mathematics and Statistics, University of Strathclyde, Glasgow G1 1XH, U.K.}


%

\begin{abstract}
In this paper, we use the truncated EM method to study the finite time strong convergence for the SDEs with
Poisson jumps under the Khasminskii-type condition. We establish the finite time $ \mathcal L ^r (r \ge 2) $ convergence  rate when the drift and diffusion coefficients satisfy super-linear condition and the jump coefficient satisfies the linear growth condition. The result shows that the optimal $\mathcal L ^r$-convergence rate is close to $ 1/ (1 + \g)$, where $\g$ is the super-linear growth constant. This is significantly different from the result on SDEs without jumps. When all the three coefficients of SDEs are allowing to grow super-linearly, the
$ \mathcal L^r (0<r<2)$ strong convergence results are also investigated and the
optimal strong convergence rate is shown to be not greater than $1/4$.
 Moreover, we prove that the truncated EM method preserve nicely the mean square exponentially stability and asymptotic boundedness of the underlying SDEs with Piosson jumps. Several examples are given to illustrate our results.
\end{abstract}

\begin{keyword}
Stochastic differential equations, local Lipschitz condition, Khasminskii-type condition, truncated EM method, Piosson jumps.
\end{keyword}

\end{frontmatter}

\linenumbers
\section{Introduction}

Due to the broad applications in modeling uncertain phenomenon, stochastic differential equations (SDEs) driven by Brownian motions have been attracting lots of attentions \cite{ALL2007a,Mao2007book}. When some unexpected events happen, some jumps may be needed to model the effects of those events. For example, a breaking news after the close of the stock market may lead to a huge difference between today's closing price and tomorrow's opening price. To take both the continuous and discontinuous random effects into consideration, SDEs driven by both Brownnian motions and Poisson jumps are often employed as a generalisation of the SDEs only driven by Brownian motions.
\par
Despite the wide applications, the explicit solutions to SDEs are hardly found. Therefore, to construct some efficient numerical methods is of extremely important. The series works of Higham and Kloeden \cite{HK2005a,HK2006a,HK2007a} studied some implicit methods for SDEs with Poisson jumps. In their papers, the strong convergence, the convergence rates and stability of different implicit methods were proposed and investigated for some SDEs, whose drift coefficient satisfies non-global Lipschitz condition, and both the diffusion coefficient and the coefficient for the Poisson jumps are global Lipschitzian.
When the global Lipschitz condition on the diffusion coefficient is disturbed, the tamed Euler and the tamed Milstein methods were proposed for SDEs driven by the more generalised process, L\'evy process \cite{DKS2016,KS2017}. The taming techniques were original proposed in \cite{Hut02} for the construction of explicit methods for SDEs with non-globally Lipschitz continuous coefficients. As indicated in \cite{Hig2011a}, explicit methods have their own advantages on the relatively simple structure and the avoidance of solving some nonlinear systems in each iteration. Therefore, the studies on explicit methods for SDEs with non-global Lipschitz coefficients have been blooming in recent years.
Sine and cosine functions were employed in \cite{ZM2017} to construct some explicit methods for SDEs with both the drift and diffusion coefficients growing super-linearly. The taming techniques were modified and generalised in \cite{Sabanis2013} and \cite{HJ2015}. The truncated Euler method were proposed in \cite{Mao2015,Mao2016}.

In this paper, we borrow the truncating idea to propose the truncated Euler method for SDEs with Poisson jumps. The main contributions of this work are twofold. Firstly, all the drift coefficient, the diffusion coefficient and the coefficient for Poisson jumps are allowed to grow super-linear. To our best knowledge, this is the first work to study an explicit numerical method for SDEs with all the three coefficients that can grow super-linearly. Secondly, both the finite time convergence and asymptotic behaviours of the method are investigated.
\par
It should be noted that the truncated Euler for SDEs with the global Lipschitzian pure jumps were studied in \cite{TY2018arXiv}. Other numerical methods for SDEs with Poisson jumps or L\'evy process were also proposed and investigated \cite{KS2017b,YW2017a,Prz2016,MYM2016,WG2010}, we just mention some of them here and refer the readers to the references therein. For the detailed and systemic introductions to numerical methods for SDEs and SDEs with jump, we refer the readers to the monographs \cite{KP1992a} and \cite{Pla2010}.
\par
This paper is constructed as follows. Section \ref{secmathpre} sees some necessary mathematical preliminaries. Section 3 contain the main results on the finite time convergence. The asymptotic behaviours, stability and boundedness, of the numerical solutions are presented in section 4. Several examples are given in the Section 5. Section \ref{secon} concludes the paper and points out some future research.

\section{Mathematical Preliminaries}\label{secmathpre}
Throughout this paper, unless otherwise specified, let $(\Omega , {\mathcal F}, \mathbb{P})$ be a complete probability
space with a filtration  $\{{\cal F}_t\}_{t\ge 0}$ satisfying the usual conditions (i.e., it is increasing and right continuous while $\cal{F}_\textrm{0} $ contains all $\mathbb{P}$\textrm{-}null sets). Let $\E$ denote the probability expectation with respect to $\mathbb{P}$. Let $B(t) $ be an $m$-dimensional Brownian motion defined on the probability space and is $ \cal{F}_\textrm{t}$-adapted. $N(t)$ is a scalar Poisson process with the compensated Poisson precess $ \widetilde N(t) =N(t) - \lambda t $, where the parameter $\lambda $ is a jump intensity. If $A$ is a vector or matrix, its transpose is denoted by $A^T$. If $x \in \R^d$, then $|x|$ is the Euclidean norm. If $A$ is a matrix, its trace norm is denoted by $ |A | = \sqrt{ (A^T A)}$. For two real numbers $a$ and $b$, we use $a \ve b = \max (a,b)$ and $a \we b = \min (a,b)$. For a set $G$, its indicator function is denoted by $\II_G$. Moreover, $ \mathcal L^r = \mathcal L^r( \Omega , {\mathcal F}, \mathbb{P})$ denotes the space of random variables $X$ with a norm $|x|_r := (\E |X|^r)^{1/r} < \infty$ for $r >0$.
  In what follows, for notational simplicity, we use the convention that $C$ represents a generic positive constant, the value of which may be different for different appearances.

Consider a $d$-dimensional SDEs with Piosson jumps:
\begin{equation}\label{eq0}
dx(t^-) = f(x(t^-))dt + g(x(t^-))dB(t) + h (x(t^-)) d N(t), \quad t \ge 0.
\end{equation}
with the initial value $x(0) = x _0 \in \R ^d $,  where $ x(t^-)$ denotes $ \lim_{s \to t^- } x(s)$. Here, $f: \R ^d \to \R^d$ is the drift coefficient, $g: \R ^d \to \R ^ {d \times m}$ is the diffusion coefficient, $h: \R ^d \to \R^d$ is the jumps coefficient.
\section{Finite time convergence}
\subsection{Convergence rate of the partially truncated EM method in $\mathcal L^r(r \ge 2)$} \label{geqthan2}
In order to discuss the convergence rates of the truncated EM method in $L^r$ for $r \ge 2$. We assume that $f$ and $g$ can be decomposed as $f(x) = F_1(x) + F(x)$ and $g(x)= G_1(x) + G(x) $, where $F_1, F: \R ^d \to \R^d$, and $G_1, G: \R ^d \to \R ^ {d \times m}$. Moreover, the coefficients $F$, $G$,  $F_1$ ,$G_1$ and $h$   satisfy the following conditions.
\begin{assumption} \label{Local_condition}
There exist positive constant $ L_1>0$ and $\g \ge 0$ such that
\begin{align} \label{Polynomial_eq}
 |F_1(x) - F_1(y) | \ve |G_1(x) - G_1(y)|  \ve |h(x) - h(y)|\le L_1|x-y| , \quad \forall  x , y \in \R ^d \no
 |F(x) - F(y)| \ve |G(x) - G(y)|  \le L_1(1+ |x|^{\g} + |y|^{\g}) |x-y|, \quad \forall  x , y \in \R ^d .
\end{align}
\end{assumption}
The parameter $\g$, which we call super-linear growth constant.
By \As \ref{Local_condition}, we can derive that there exists a positive constant $K_1$ such that
\begin{align}\label{linear_eq}
|F_1(x)|\ve  |G_1(x)| \ve |h(x)|  \le K_1 (1+ |x|), \qu \forall x \in \R ^d,
\end{align}
 which implies that $F_1$, $G_1$  and $h$ satisfy the linear growth condition. Similarly, we have
 \begin{align}\label{fgh}
|F(x)| \ve |G(x) | \le (2L_1+ |F(0)|+ |G(0)|)|x| ^ {1+ \g} , \qu \forall |x | \ge 1.
 \end{align}
We also impose the following standing hypotheses.
\begin{assumption} \label{onesided_conditions}
 There exists a pair of  constant $ \bar r > 2 $ and $L_2 >0$  such that
\begin{align} \label{onesided_eq}
( x - y )^T( F(x) -F(y))   + \frac{\bar r-1}{2} |G(x) - G(y)|^2 \le L_2 |x -y | ^2 ,\quad \forall \; x,y \in \R ^d.
\end{align}
\end{assumption}
By \As \ref{onesided_conditions}, we can derive that for any $ r \in [2, \bar r)$
\begin{align}\label{eq29}
 (x-y)^T (f(x) -f(y))  + \frac{r-1}{2} |g(x)-g(y)|^2 \le L_3 |x-y|^2.
\end{align}
where $L_3 = 2L_1 + L_2 + \frac{L_1^2 + (r-1)(\bar r -1 )}{ \bar r -r}$ (see \cite{Guo2017partial}).
\begin{assumption} \label{monotone_condition}
(Khasminskii-type condition) There exist constants $\bar p > \bar r  $, $K_2>0$ such that
\begin{align}\label{monto_eq}
 x ^T  F(x)  + \frac{\bar p-1}{2} | G(x)|  ^2  \le K_2(1 + |x|^2 ), \quad  \forall x \in \R ^d.
\end{align}
\end{assumption}
By \As \ref{monotone_condition}, we also have that for any $ p \in [2, \bar p) $
\begin{align}\label{eq28}
 x ^T  f(x)  + \frac{p-1}{2} |g(x)|^2 \le K_3 (1+ |x|^2),
\end{align}
where $K_3 = 2K_1 + K_2 + \frac{K_1^2 + (p-1)(\bar p -1 )}{ \bar p -p}$ (see \cite{Guo2017partial}).

The truncated idea is to deal with the super-linear coefficients. In the viewpoint of the finite-time convergence, the linear coefficient does not cause any problem to the EM scheme and hence there is no need to truncate it \cite{Guo2017partial}. In our truncated EM method, we only truncate the super-linear terms, that is $F$ and $G$.
To define the truncated EM scheme, we first choose a strictly increasing function $\mu : \R ^+  \rightarrow \R^ +$ such that $\mu (n) \rightarrow \infty $, as $ n \rightarrow \infty$,
and
\begin{align*}
\sup _{|x| \le n}  |F(x)| \vee |G(x)| \le \mu (n), \quad \forall n \ge 1.   
\end{align*}
Denoted by $\mu ^{-1}$ is the inverse function of $\mu$. We also choose  a strictly decreasing function $ \ph: (0, 1) \rightarrow (0, \infty)$ such that
 \begin{align}\label{delta_h_relations}
\lim _ {\Delta \to 0} \ph(\Delta) = \infty \quad  \textrm{  and} \quad  (\ph(\DD))^{\bar p} \le  \DD ^ {-1} \we \DD ^ {-\bar p /4},   \qu \forall \DD \in (0,1] .
 \end{align}
For a given step size $\Delta \in (0, 1)$, let us define a mapping $\pi _{\Delta }$ from $\R ^d$ to the closed ball $\{ x \in \R ^d : |x| \le \mu^{-1} (\ph(\Delta )) \}$ by
\begin{align*}
 \pi _{\Delta} = \big (  |x|  \wedge \mu^{-1} (\ph(\Delta )) \big  ) \frac{x}{|x|} .
       \end{align*}
We set $x / |x| = 0$ when $x=0$. We then define the partially truncated functions
\begin{align*}
F_{\Delta} (x)= F( \pi _ {\Delta}(x)) , \quad
 G_{\Delta}(x)  = G( \pi _ {\Delta}(x)) , \quad \forall x \in \R^d \\
 \fD (x) = F_1(x)+F_{\Delta} (x) \quad \textrm{and} \quad  \gD (x) = G_1(x)+F_{\Delta} (x), \quad \forall x \in \R^d
\end{align*}
It is easy to see that
 \begin{align} \label{HD}
|F_{\Delta}(x)| \vee |G_{\Delta}(x)|  \le \ph(\Delta), \quad \forall x \in \R ^d.
\end{align}
Obviously, $F_{\Delta}$ and $G_ {\Delta}$ are bounded while $F$ and $G$ may not. The following lemma shows that the truncated functions maintain the Khaminskii-type condition nicely (see \cite{Mao2015}).
\begin{lemma} \label{lem 24}
Let Assumption \ref{monotone_condition} hold. Then, for all $\Delta \in (0, 1]$, we have
\begin{align*}
x^{T} F_ {\Delta} (x) + \frac{\bar p-1}{2} |G_{\Delta}(x)|^2 \le 2K_2 (1+ |x|^2), \quad \forall x \in \R^d.
\end{align*}
\end{lemma}
We can show that for any $ p \in [2, \bar p)$ , we have
\begin{align}\label{eq209}
x^{T} f_ {\Delta} (x) + \frac{p-1}{2} |g_{\Delta}(x)|^2 \le K_4 (1+ |x|^2), \quad \forall x \in \R^d.
\end{align}
where $K_4 = 2K_1 + 2K_2 + \frac{K_1^2 + (p-1)(\bar p -1 )}{ \bar p -p}$ (see \cite{Guo2017partial}).
We now form the discrete-time truncated EM numerical solutions $X_{\Delta}(t_k^-)  \approx x(t_k^-)$, for $t_k^- = k \Delta t$ by setting $X_ {\Delta} (0) = x_0$ and computing
\begin{align}\label{EM scheme}
\XD (t_{k+1}^- ) = \XD (t_k^-) + \fD(\XD(t_k^-)) \DD t + \gD (\XD (t_k^-)) \DD B_k + h (\XD (t_k^-)) \DD N_k, \quad 0 \le k \le M-1,
\end{align}
where $\DD B_k = B(\tkk) - B (\tk)$, $\DD N_k = N(\tkk) - N (\tk)$. It is consentient to use the continuous-time step process $\bar x_{\Delta} (t^-) $ which is defined by
\begin{align}\label{step_process}
\bxdt = \sum _ {k=0}^ {\infty} \XD (\tk) \II_{[\tk,\tkk)} (t)
\end{align}
where $\II$ is a indicator function. The other continuous-time process is defined by
\begin{align*}
\xdt  = x_0 + \int_0^t \fdxs ds + \int_0^t \gdxs dB(s) + \int_0^t h(\bxds) dN(s).
\end{align*}
It is easy to see that $x_{\Delta}(\tk) = \bar x _ {\Delta} (\tk)  = \XD (\tk)$. Moreover, $\xdt  $ is an It\^{o} process with It\^{o} differential
\begin{align*}
d \xdt = \fdxt dt + \gdxt dB(t) + h(\bxdt) dN(t).
\end{align*}
We first state a known result (see \cite{DKS2016}) as a lemma.
\begin{lemma} \label{lemma4.77}
Under Assumption \ref{Local_condition} and  \ref{monotone_condition} the SDE \eqref{eq0} has a unique global solution $x(t)$, moreover, for any $p \in [2, \bar p)$,
\begin{align*}
\sup _ {0 \le t \le T} \E |x(t)| ^ { p} < \infty , \quad \forall T > 0.
\end{align*}
\end{lemma}
In order to bound the $p$-th moment of the truncated EM solution, we need the following lemma.
\begin{lemma} \label{lem17}
For any $\DD \in (0, 1]$ and $t >0$, we have
\begin{align}\label{eq114}
\E | \xdt  - \bxdt|^{\hat p} &  \le C \Big ((\ph(\DD))^{\hat p}\DD ^ {\hat p  /2 }+ (1+ \E|\bxdt|^{\hat p})\DD \Big ), \quad  \forall {\hat p} \ge 2.
                                                \end{align}
\end{lemma}
\pr  Fix any $\DD \in (0, 1]$, $t \ge 0$ and $\hat p \ge 2$. There is a integer $k \ge 0$ such that $\tk \le t < \tkk$. By \As \ref{Local_condition} and \eqref{HD}, we have
\begin{align}\label{eq116}
& \E | \xdt - \bxdt|^{\hat p}  \\ \nonumber
 & = \E \left | \int_{\tk}^t \fdxs ds + \int_{\tk}^t \gdxs dB(s) +\int_{\tk}^t  h(\bxds) dN(s) \right |^{\hat p} \\ \nonumber
 & \le C \left ( \E \left | \int_{\tk}^t \fdxs ds \right |^{\hat p} + \E \left | \int_{\tk}^t \gdxs dB(s) \right|^{\hat p}
    + \E \left | \int_{\tk}^t h(\bxds) dN(s) \right |^{\hat p}      \right )\\ \nonumber
 & \le C \Big ( \DD ^ {\hat p-1} \E  \int_{\tk}^t\left | \fdxs  \right |^{\hat p} ds + \DD ^ {(\hat p-2)/2 }\E  \int_{\tk}^t \left | \gdxs  \right|^{\hat p} ds +
                      \E \left | \int_{\tk}^t h(\bxds) d N (s) \right |^{\hat p}   \Big ) \\ \nonumber
 &  \le C \Big (\DD ^ {\hat p/2}( 1+ \E | \bxdt|^{\hat p} + (\ph(\DD))^{\hat p} )    +
                      \E \left | \int_{\tk}^t h(\bxds) d N (s) \right |^{\hat p}   \Big ),
\end{align}
where $C$ is a generic constant, the value of which may change between occurrences.
By the characteristic function's argument \cite{Bao2011}, for $\DD \in (0,1]$
\begin{gather}\label{memont of Nk}
\E |\DD N_k|^{\hat p} \le c_0 \DD,
\end{gather}
where $c_0$ is a positive constant which is independent of $\DD$.
Therefore,
\begin{align*}
 & \E \left | \int_{\tk}^th(\bxds) d N (s) \right |^{\hat p}   = \E | h (\XD (t_k^-)) \DD N_k|^ {\hat p} \no
 & =  \E | h (\XD (t_k^-))| ^ {
\hat p}  \E | \DD N_k|^ {\hat p}
                              \le C (1+ \E|\bxdt|^{\hat p})\DD.
        \end{align*}
Inserting this into \eqref{eq116} and combing with $\DD^{\hat p / 2} \le \DD $ gives
\begin{align*}
\E | \xdt - \bxdt |^{\hat p} \le C \Big ((\ph(\DD))^{\hat p}\DD ^ {\hat p  /2 }+ (1+ \E|\bxdt|^{\hat p})\DD \Big ).
\end{align*}
Thus, we complete the proof. $\Box$

\begin{lemma} \label{moment of xdt}
Let \As \ref{Local_condition} and \ref{monotone_condition} hold and let $ p \in [2, \bar p) $ be arbitrary. Then
\begin{align}\label{eq moment lof xdt }
\sup _{ 0 \le \DD \le 1} \sup _ {0\le t \le T}  \E | \xdt|^p \le C, \qu \forall T > 0,                                                                        \end{align}
\end{lemma}

\pr Fix any $ \DD \in (0,1]$ and $T >0 $. By the It\^{o} formula and  \eqref{eq209}, we have
\begin{align}\label{eq119}
  \E | \xdt|^p - |x_0|^p & \le    \E \int_0^t p |\xdt|^{p-2}\Big (  x^T _ {\DD} (s) \fdxs  + \frac{p-1}{2} |\gdxs|^2  \Big )ds \no
   & \qu + \lambda \E \Big ( \int_0^t | \xds + h(\bxds) |^p - | \xds |^p  \Big )ds \no
    & \le  \E \int_0^t p |\xdt|^{p-2}\Big (  \bar x^T _ {\DD} (s) \fdxs  + \frac{p-1}{2} |\gdxs|^2  \Big )ds \no
    & \qu +   \E \int_0^t p |\xds|^{p-2} (\xds - \bxds )^T \fdxs ds  \no
    & \qu + \lambda \E \Big ( \int_0^t | \xds + h(\bxds) |^p - | \xds |^p  \Big )ds  \no
    & \qu \le I_1 + I_2 + I_3+ I_4,
\end{align}
where \begin{align}\label{eq306}
I_1 & = \E \int_0 ^ t p  K_4 |\xds|^{p-2} (1+ |\bxds|^2) ds , \\
I_2  & = \E \int_0^t p |\xds|^{p-2} |\xds - \bxds||F_1 (\bxds)| ds ,\\
I_3  & = \E \int_0^t p |\xds|^{p-2} |\xds - \bxds||F_{\DD} (\bxds)| ds ,
      \end{align}
and
\begin{align}\label{eq307}
I_4  & = \lambda \E \Big ( \int_0^t | \xds + h(\bxds) |^p - | \xds |^p  \Big )ds  .
\end{align}
By the Young inequality
\begin{align*}
a ^ {p-2} b^2 \le \frac{p-2}{p} a ^p + \frac{2}{p} b ^ {p}, \qu \forall a,b \ge 0,
\end{align*}
we then have
\begin{align}\label{310}
I_1 \le C \Big ( 1+ \intt ( \E |\xds|^p + \E |\bxds|^p )ds  \Big ).
\end{align}
Similarly, we can show that
\begin{align}\label{311}
I_2 \le C \Big ( 1+ \intt ( \E |\xds|^p + \E |\bxds|^p )ds  \Big ).
\end{align}
By \As \ref{Local_condition}, it is not difficult to prove that there exists a positive constant $c_1$ such that \begin{align}\label{eq120}
| \xds + h(\bxds) |^p - | \xds |^p \le c_1(1+  |\xdt|^p + |\bxdt |^p).
\end{align}
Hence, we have
\begin{align}\label{312}
I_4 \le C \Big ( 1+ \intt ( \E |\xds|^p + \E |\bxds|^p )ds  \Big ).
\end{align}
Moreover, by the Young inequality and \eqref{HD}, we get
\begin{align}\label{eq121}
I_3 \le (p-2) \E \intt |\xds|^p +  2  \E \intt  |\xds - \bxds|^{p/2} |F_{\DD} (\bxds )|^{p/2} ds.
 \end{align}
By Lemma \ref{lem17} and \eqref{HD}, we obtain
\begin{align}\label{eq1222}
& \E \int_0^T |\xds -\bxds |^{p/2} | F_{\DD} (\bxds )|^{p/2}ds  \no
& \le (\ph(\DD))^{p/2} \int_0^T \E | \xds - \bxds|^{p/2}ds  \no
& \le (\ph(\DD))^{p/2} \int_0^T (\E | \xds - \bxds|^{p}) ^ {\frac{1}{2}}ds  \no
& \le C (\ph(\DD))^{p/2} \Big (  (1 + \E | \bxds|^{p/2})\DD^ {1/2} +  (\ph(\DD))^{p/2} \DD ^ {p/4} \Big ) \no
& \le C  \Big (  (1 + \E | \bxds|^{p})(\ph(\DD))^{p/2} \DD^ {1/2}+  (\ph(\DD))^{p} \DD ^ {p/4} \Big ) \no
& \le C (1 + \E | \bxds|^{p}).
\end{align}
Noting that the last inequality in \eqref{eq1222} has used the condition \ref{delta_h_relations} which implies
\begin{align*}
(\ph(\DD))^{p/2} \DD^ {1/2} \le 1 \qu \textrm{and} \qu  (\ph(\DD))^{p} \DD ^ {p/4} \le 1.
\end{align*}
Inserting \eqref{eq1222} into \eqref{eq121} also gives
\begin{align}\label{313}
I_3 \le C \Big (1+  \intt ( \E |\xds|^p + \E |\bxds|^p )ds  \Big ).
\end{align}
Substituting \eqref{310}, \eqref{311}, \eqref{312} and \eqref{313} into \eqref{eq119}, we get
\begin{align*}
 \E |\xdt|^p  & \le C \Big ( \intt ( 1+ \E |\xds|^p + \E |\bxds|^p )ds  \Big ) \no
  & \le C  \Big ( 1 +  \intt   \sup _ {0 \le u \le s}  \E | x _ {\DD} (u)|^p ds \Big ) .
\end{align*}
Then, we have
\begin{align*}
 \sup _ {0 \le u \le t }\E |x_{\DD}(u^-)|^p \le C  \Big ( 1 +  \intt   \sup _ {0 \le u \le s}  \E | x _ {\DD} (u)|^p ds \Big ) .
\end{align*}
The Gronwall inequality yields
$$  \sup _ {0 \le u \le T} \E  |  x _ {\DD} (u^-)|^p \le C.$$
As this holds for any $\DD \in (0, 1]$ while $C$ is independent of $\DD$, we obtain the required assertion. $\Box$

The following lemma shows that $x_{\DD}(t)$ and $\bar x _ {\DD}(t)$ are close to each other in the sense of $\mathcal L ^p$.
\begin{lemma} \label{lemma33}
Let \As \ref{Local_condition} and \ref{monotone_condition} hold and let $ t \in [0, T] $. Then there is a $\bar \DD   \in (0, 1]$ such that for all $ \DD \in (0, \bar \DD]$,
\begin{align}\label{temp2}
\E | \xdt  - \bxdt|^{ p} &  \le C \B ( (\ph(\DD))^{ p}\DD ^ { p  /2 } + \DD \B), \quad  2 \le p < \bar p, \\
\E | \xdt  - \bxdt|^{ p} &  \le C \B ( (\ph(\DD))^{ p } \DD ^ {p /2}   + \DD ^ { p/2}\B ),\quad  0 < p <2 .
\end{align}
Consequently
\begin{align}\label{eq115}
\lim _ {\DD \to 0} \E |\xdt - \bxdt|^{ p} =0, \qu p >0.
\end{align}
\end{lemma}
\pr   For any $p \ge 2$, by Lemma \ref{moment of xdt}, there is a $\bar \DD \in (0, 1]$ such that
\begin{align} \label{temp1}
\sup _{ 0 \le \DD \le \bar \DD} \sup _ {0\le t \le T}  \E | \xdt|^p \le C.                                                                        \end{align}
Now, fix any $\DD \in (0, \bar \DD ]$, inserting \eqref{temp1} into \eqref{eq114} gives \eqref{temp2}. For any $ p \in (0,2)$, the \Holder implies
\begin{align*}\label{}
& \E | \xdt - \bxdt |^{ p} \le \Big ( \E | \xdt - \bxdt |^2  \Big ) ^ { p /2} \no
& \le  C \Big(   (\ph(\DD))^{2} \DD  + \DD \Big ) ^ { p /2}  =  C \B ( (\ph(\DD))^{ p } \DD ^ {p /2}   + \DD ^ { p/2}\B ).
\end{align*}
Noting from \eqref{delta_h_relations} that $ (\ph ( \DD))^{ p}\DD ^ { p  /2 } \le \DD ^ {p/4} \we \DD ^ {p/2-1}$ for $ 2 < p <\bar p$, we obtain \eqref{eq115} from \eqref{temp2}. $\Box$

Let us propose two lemmas before we state our main results in this paper.
\begin{lemma} \label{lemma303}
Let \As \ref{Local_condition} and \ref{monotone_condition} hold. For any real number $ n  > |x_0| $, define the
stopping time \begin{align*}
 \tr= \inf \{ t \ge 0 : |x(t)| \ge n \}.
              \end{align*}
              Then
              \begin{align}\label{eq35}
\PP (\tr \le T)  \le  \frac{C}{n^2}.
              \end{align}
\end{lemma}
\pr The proof is given in the Appendix. $\Box$
\begin{lemma} \label{lemma34}
Let \As \ref{Local_condition} and \ref{monotone_condition} hold. For any real number $n > |x_0| $, define the
stopping time,
\begin{align*}
 \rho _ {\DD , n}= \inf \{ t \ge 0 : |\xdt| \ge n \}.
              \end{align*}
              Then
              \begin{align}\label{eq35}
\PP (\rho _{\DD , n } \le T)  \le  \frac{C}{n^2}.
              \end{align}
\end{lemma}
\pr The proof is given in the Appendix. $\Box$

Now, we will show one of our main results in our paper. The proof is similar to that of Theorem 3.6 in \cite{Guo2018note}, we only highlight the different parts.
\begin{theorem} \label{lemma3.2}
Let \As \ref{Local_condition}, \ref{onesided_conditions} and \ref{monotone_condition} hold and assume that there exists a number $ p \in (2, \bar p)$ such that
\begin{align}\label{ass2}
 p > (1+ \g )\bar r  .
\end{align}
Let $ r \in [2, \bar r)$ be arbitrary. Then for any $\DD \in (0,1]$,
\begin{align}\label{resu1}
\E | x(T) - x_{\DD}(T)|^r \le C \B ( (\mu ^{-1} (\ph(\DD)))^{-(p- (1+ \g)r)}  +  (\ph(\DD))^r \DD ^ {r/2} + \DD ^ {(p - \g r )/p}\B )
\end{align}
and
\begin{align}\label{resu2}
\E | x(T) - \bar x_{\DD}(T)|^r \le C \B ( (\mu ^{-1} (\ph(\DD)))^{-(p- (1+ \g)r)}  +  (\ph(\DD))^r \DD ^ {r/2} + \DD ^ {(p - \g r )/p} \B ).
\end{align}
In particular, we define
\begin{align}\label{mu_u}
\mu (x) = L_4 x ^{1+ \g}, \qu x \ge 0,
\end{align}
where $L_4 = 2L_1 + |F(0)|+ |G(0)|$, and let
\begin{align}\label{mu_2}
\ph(\DD)  =  \DD ^ {- \veps} \qu \textrm{for } \; \textrm{some} \qu  \veps  \in (0, 1/4 \we  1/p ]
\end{align}
to obtain
\begin{align}\label{resu3}
\E | x(T) - x_{\DD}(T)|^r \le C  \DD ^ { [\veps (p - (1+ \g ) r)/ (1+ \g)] \we [r (1-2 \veps)/2 ]  \we [(p - \g r )/p]}
\end{align}
and
\begin{align}\label{resu4}
\E | x(T) - \bar x_{\DD}(T)|^r \le C  \DD ^ { [\veps (p - (1+ \g ) r)/ (1+ \g)] \we [r (1-2 \veps)/2 ] \we (p - \g r )/p }
\end{align}
for all $ \DD \in (0,1]$.
\end{theorem}
\pr Let $\DD \in (0,1]$ be arbitrary.  Let $ e_{\DD}(t) = x(t) - x_{\DD}(t)$ for $t >0$. Fix a number $q \in (r ,\bar r)$, \eqref{ass2} means $ p > (1 + \g) q $.
For any integer $n > |x_0|$, define the stopping time
\begin{align*}
\sigma _n  = \inf \{  t \ge 0 : |x(t)| \ve |x_ {\DD}(t) |  \ge n \}.
\end{align*}
By the It\^{o} formula, we get that for $ 0 \le t \le T$
\begin{align}\label{eq43}
& \E |e _ {\DD} (t \we \sigma _n )|^r   \no
& \le \E \int_0^{t \we \sigma _n } r | \eds |^{r-2} \Big ( e^T _ {\DD}(s^-) (f(x(s^-) )- \fdxs) + \frac{r-1}{2} |g(x(s^-))  -\gdxs | ^2 \Big ) ds\no
&  \qu + \la \E \int_0^{t \we \sigma _n }    \Big ( |\eds + (h(x(s^-))- h(\bar x _ {\DD}(s^-))) |^r  - | \eds |^r \Big )ds \no
& =: J_1 + J_2.
\end{align}
Let us estimate $J_2$ first. Using \As \ref{Local_condition} gives
\begin{align*}
 & | x(s^-) - \xds  + h(x(s^-)) - h(\bxds ) |^r \no
& \le 2^{r-1}(|x(s^-) - \xds| ^r  + |h(x(s^-)) - h(\bxds )|^r)  \no
& \le 2^{r-1}(|x(s^-) - \xds| ^r  + L_1 ^r |x(s^-) - \bxds |^r)  \no
& \le c_1 ( | x(s^-) - \xds |^r  + |\xds -\bxds |^r  ),
\end{align*}
where $c_2 = 2^{r-1} (1+ L_1^r 2 ^ {r-1}) >1$.
Hence, by Lemma \ref{lemma33}, we have
\begin{align}\label{eq 50}
J_2 & \le \la (c_2 -1)  \intt \E | e_{\DD} (s \we \sigma _n  )|^r ds + \la c_2  \int _0 ^ T \E | \xds - \bxds |^r ds \no
   & \le \la (c_2 -1) \intt \E | e_{\DD} (s \we \sigma _n )|^r ds + C \B ( \DD ^ {r/2} (\ph(\DD))^ {r} + \DD \B ).
\end{align}
By the elementary inequalty, $J_1$ can be decomposed into two parts denoted by $J_1 = J_3 + J_4$, where
\begin{align}\label{J3}
J_3 =  \E \int_0^{t \we \sigma _n } r | \eds |^{r-2}  & \Big ( e^T _ {\DD}(s^-) (f(x(s^-) )- f(\xds)) \no
    & + \frac{q-1}{2} |g(x(s^-))  -g(\xds) | ^2 \Big ) ds\
\end{align}
and
\begin{align}\label{J4}
J_4 =  \E \int_0^{t \we \sigma _n } r | \eds |^{r-2}  & \Big ( e^T _ {\DD}(s^-) (f(\xds)  )- f_{\DD}( \bar x _ {\DD}(s ^-)  )) \no
    & + \frac{(r-1)(q-1)}{2(q -r)} |g(\xds)  - g_{\DD}( \bar x _ {\DD}(s ^-)  )| ^2 \Big ) ds.
\end{align}
By \eqref{eq29}, we have
\begin{align}\label{JJ1}
J_3 \le  r L_3  \int_0^{t \we \sigma _n }  \E | \eds |^{r}  ds.
\end{align}
The Elementary inequality gives
\begin{align}\label{JJ2}
J_4 &  \le   \E \int_0^{t \we \sigma _n } r | \eds |^{r-2}   \Big ( e^T _ {\DD}(s^-) (f(\xds)  )- f_{\DD}(\xds)) \no
    &  \qu + \frac{(r-1)(q-1)}{(q -r)} |g(\xds)  - g_{\DD}(\xds)| ^2 \Big ) ds  \no
 & +    \E \int_0^{t \we \sigma _n } r | \eds |^{r-2}   \Big ( e^T _ {\DD}(s^-) (f_{\DD}(\xds)  )- f_{\DD}(  \bar x_ {\DD}(s ^-)  )) \no
    & \qu  + \frac{(r-1)(q-1)}{(q -r)} |g_{\DD}(\xds)  - g_{\DD}( \bxds)| ^2 \Big ) ds \no
 &   = : J_{41} + J_ {42}.
\end{align}
In the same way as Theorem 3.6 in \cite{Guo2018note} was proved, we can show that
\begin{align}\label{resu51}
J_{41} &  \le    C   \B (  \int _0 ^{ t \we \sigma _n } \E |e_ {\DD}(s^-)|^r ds +     (\mu ^{-1} (\ph(\DD)))^{-(p- (1+ \g)r)}  \B )
\end{align}
and
\begin{align}\label{resu52}
J_{42} &  \le    C     \int _0 ^{ t \we \sigma _n } \E |e_ {\DD}(s^-)|^r ds +  C \int _ 0 ^ T  \B(   \E |\xds - \bxds |^{ pr /( p - \g r)}  \B) ^ { (p - \g r)/ p} \no
&  \le    C     \int _0 ^{ t \we \sigma _n } \E |e_ {\DD}(s^-)|^r ds + C \B ( (  \ph (\DD) )^ { pr / ( p - \g r)} \DD ^ { 0.5pr / ( p - \g r)}  + \DD  \B )^ { ( p - \g r) /p}  \no
&  \le    C     \int _0 ^{ t \we \sigma _n } \E |e_ {\DD}(s^-)|^r ds + C \B ( (  \ph (\DD) )^ { r} \DD ^ { r/2}  + \DD ^ {( p - \g r) /p }  \B )  ,
\end{align}
where we use the Lemma \ref{lemma33} and the fact that $$ \frac{ p r}{ p - \g r} = r \frac{p}{ p - \g r} >2.$$
Inserting \eqref{resu51} and \eqref{resu52} into \eqref{JJ2}, we have
\begin{align}\label{JJ5}
J_4  \le  C   \B (  \int _0 ^{ t  } \E |e_ {\DD}(s \we \sigma _n )|^r ds +     (\mu ^{-1} (\ph(\DD)))^{-(p- (1+ \g)r)}  + (\ph(\DD))^r \DD ^ {r/2}   + \DD ^ {( p - \g r) /p }  \B ).
\end{align}
Combing \eqref{eq 50}, \eqref{JJ1} and \eqref{JJ5}, we have
\begin{align*}
\E |e _ {\DD} (t \we \sigma _n )|^r &   \le  C   \B (  \int _0 ^{ t  } \E |e_ {\DD}(s \we \sigma _n )|^r ds +     (\mu ^{-1} (\ph(\DD)))^{-(p- (1+ \g)r)}   \no
  &\qu  + (\ph(\DD))^r \DD ^ {r/2}   + \DD ^ {( p - \g r)/p}  \B ).
\end{align*}
The Gronwall inequality implies
\begin{align*}
\E |e _ {\DD} (T \we \sigma _n )|^r   \le  C   \B (   (\mu ^{-1} (\ph(\DD)))^{-(p- (1+ \g)r)}  + (\ph(\DD))^r\DD ^ {r/2}  + \DD ^ {( p - \g r)/p} \B ).
\end{align*}
Using Lemma \ref{lemma303} and \ref{lemma34} and letting $n \to \infty$ gives the desired assertion \eqref{resu1}.  By \eqref{resu1} and Lemma \ref{lemma33} gives the another assertion \eqref{resu2}. Recalling \eqref{mu_u}, then $ \mu ^{-1} = (x / L_4)^ {1/ (1 + \g)}$. Substituting this and \eqref{mu_2} into \eqref{resu1} gives \eqref{resu3}. Similarly, we can get \eqref{resu4}. Thus, the proof is complete. $\Box$
\begin{corollary}\label{cor1}
Let \As \ref{Local_condition}, \ref{onesided_conditions} hold and let \As \ref{monotone_condition} holds for all $\bar p \in ( \bar r, \infty)$.  Let $\mu$ and $\ph$ be defined in \eqref{mu_u} and \eqref{mu_2}. Then, for any
\begin{align}\label{conddd}
r \in [2,\bar r), \qu p \in ((1+ \g ) r \ve \bar r ,\bar p ) \qu \textrm{and} \qu  \veps \in (0, 1/4 \we 1/p],
\end{align}
we have
\begin{align}\label{resu6}
\E | x(T) - x_{\DD}(T)|^r \le C  \DD ^ { \veps (p - (1+ \g ) r)/ (1+ \g)  \we ( p - \g r)/p  }
\end{align}
and
\begin{align}\label{resu7}
\E | x(T) - \bar x_{\DD}(T)|^r \le C  \DD ^ { \veps (p - (1+ \g ) r)/ (1+ \g)  \we ( p - \g r)/p  }.
\end{align}
\end{corollary}
\pr Condition \ref{conddd} means $$ \veps \le \frac{1}{p} < \frac{r(1+ \g)}{ 2 p} $$   which lead to  
$$\veps (p - (1+ \g ) r)/ (1+ \g)  < r ( 1- 2 \veps)/2 . $$
By \eqref{resu3} and \eqref{resu4}, we obtain \eqref{resu6} and \eqref{resu7}. $\Box$
\begin{remark}
Replacing condition \eqref{ass2}, that is $ p > (1+ \g )\bar r $, by a weaker one $ p > (1+ \g ) r  \ve \bar r $ does not
affect the results in Theorem \ref{lemma3.2}. But, this small change will make the choice of $p$ more flexible in practice.
\end{remark}

\begin{figure}[!t]
\centering
\begin{minipage}[h]{0.45\textwidth}
   \centering
   \includegraphics[width=6cm,height=5cm]{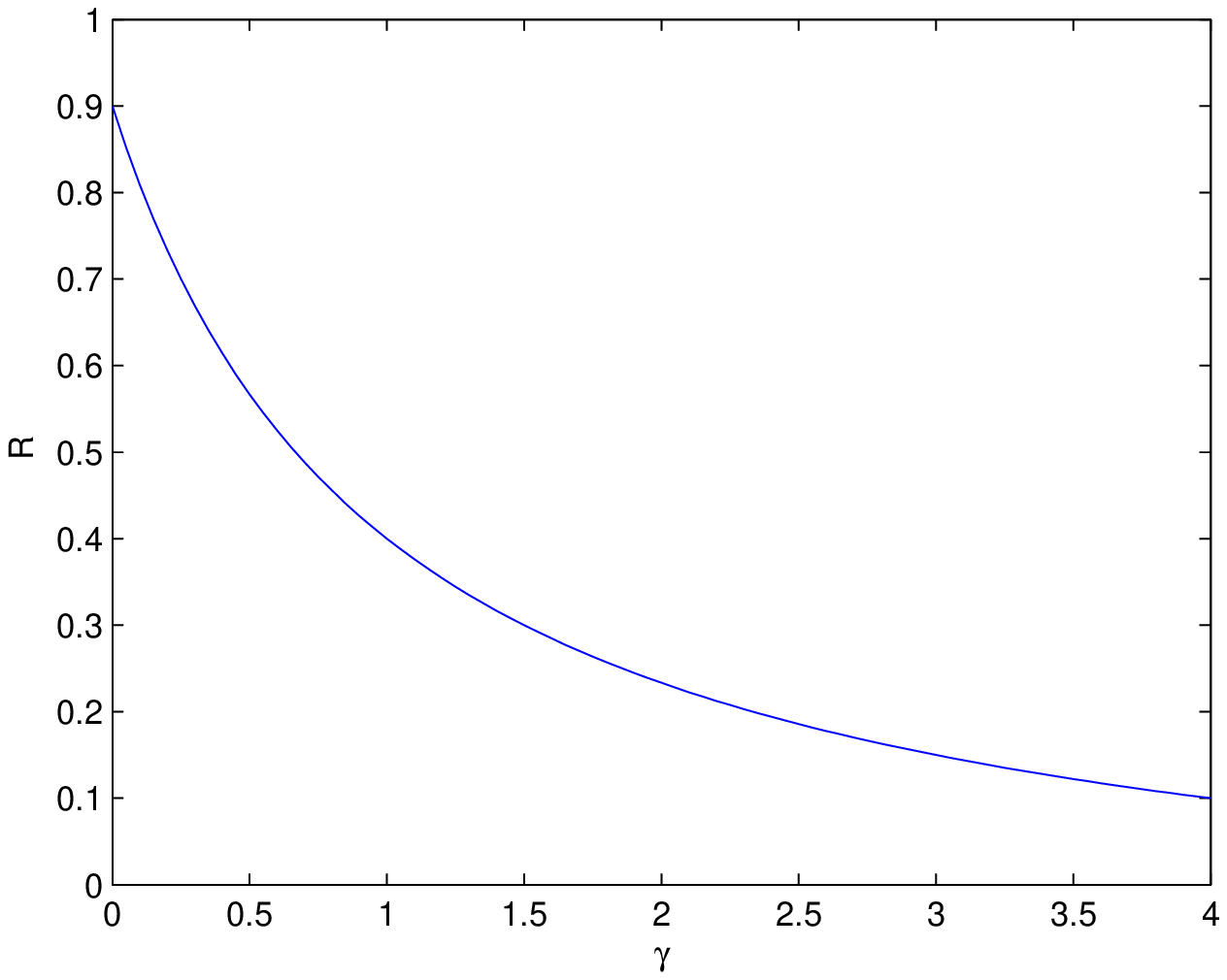}
\caption{Convergence rate $R$ versus growth constant $\g$ with $p=20$, $\veps = 1/p$ and $r=2$ }
\label{fig1}
\end{minipage}
\hfil
\begin{minipage}[h]{0.45\textwidth}
\centering
  \includegraphics[width=6cm,height=5cm]{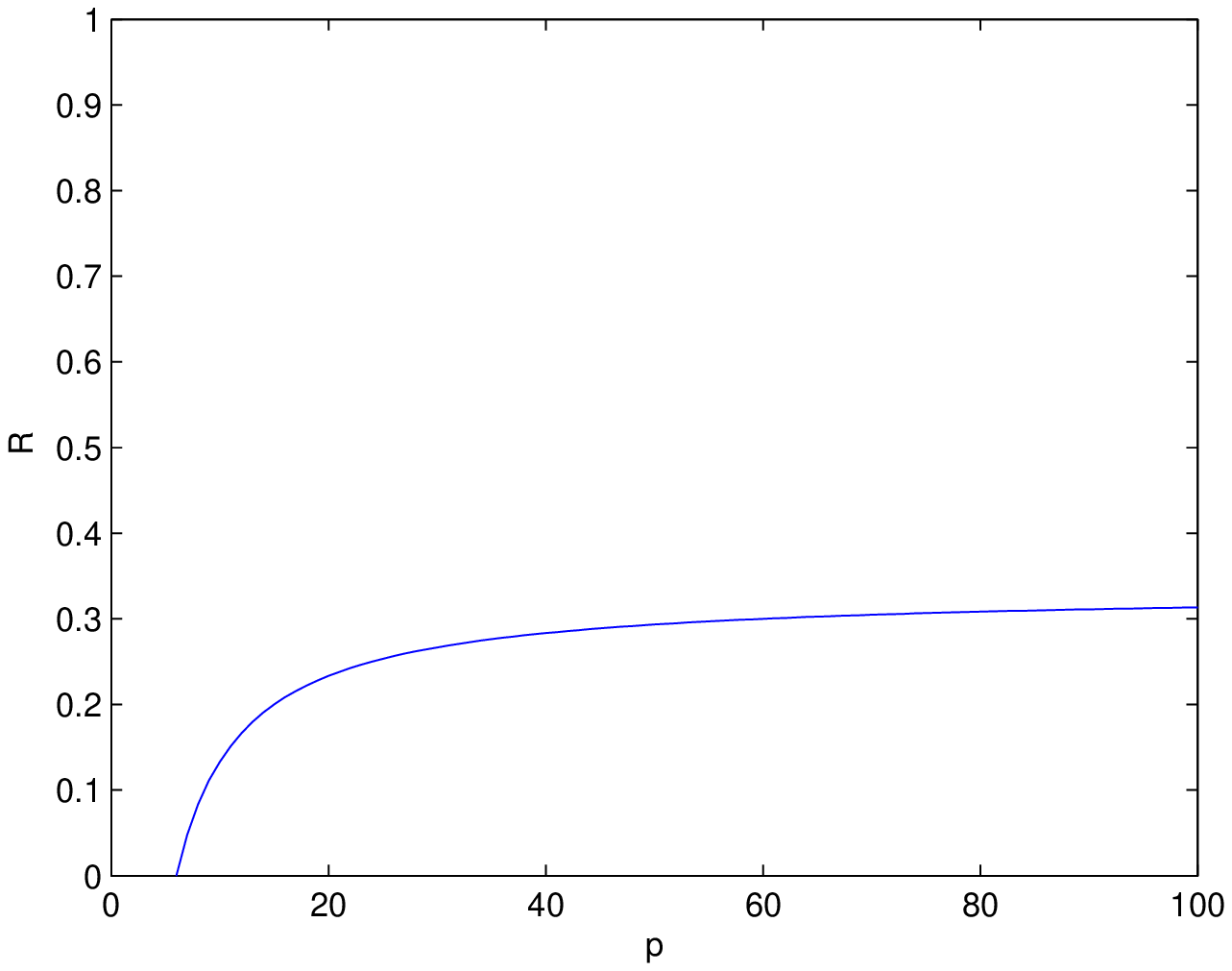}
\caption{Convergence rate $R$ versus $p$  with $\g =2 $, $\veps = 1/p$ and $r=2$}
\label{fig2}
\end{minipage}
\end{figure}

\begin{remark}
Without loss of generality, we assume that $ \bar p $ is allowed to sufficiently large. In the following discussion, we fixing $2 \le r < \bar r$.
By Corollary \ref{cor1}, we can conclude that when $p$, which is a parameter in the Khaminskii-type condition defined in \eqref{eq28}, is sufficiently large relative to $\g r$, then the order of $\mathcal L ^r$-convergence of the truncated EM method is mainly determined by the expression which we call the convergence rate function, namely
\begin{align}\label{rate_function}
R(p, \g) : =\frac{\veps (p - (1+ \g)r)}{ 1+ \g},  \qu  \forall  p \in ((1 + \g ) r \ve \bar r, \bar p),\qu \forall \g \in [0, \infty),
\end{align}
for some $ 0 < \veps \le 1/4 \we 1/p$. Noting that $R $ is proportional linearly to $ \veps$, which means that letting $  \veps = 1/4 \we 1/p$ will make convergence rate $R$ as large as possible. Hence, this setting gives
\begin{align}\label{rate_func2}
R(p, \g)  =\frac{(1/4 \we 1/p) (p - (1+ \g)r)}{ 1+ \g}.
\end{align}
To discuss the optimal convergence rate function \ref{rate_func2}, let us consider the following two possible cases:
\begin{itemize}
  \item If $ 2 < p <4 $ and $ ( 1+ \g ) r  \ve \bar r < p < \bar p$, then \eqref{rate_func2} gives
  \begin{align}\label{rate11}
   R(p, \g)  = \frac{1}{4} \B ( \frac{p}{ 1+ \g} - r \B) = \frac{p}{4} \frac{1}{ 1+ \g} - \frac{r}{4} < \frac{1}{1+ \g}.
  \end{align}
  \item If $ p \ge 4$ and $ (1+ \g ) r \ve \bar r < p  < \bar p$, then \eqref{rate_func2} gives
  \begin{align}\label{rate22}
  R(p, \g) = \frac{ p - ( 1 + \g )r}{p ( 1 + \g)} = \frac{1}{ 1+ \g} - \frac{r}{p} < \frac{1}{ 1+ \g}.
  \end{align}
\end{itemize}
From \eqref{rate22}, we can see that $R(p,\g)$ is a decreasing function with respect to $\g$ and is a increasing function with respect to $p$.  Fig \ref{fig1} and Fig \ref{fig2} show this relationship between convergence rate and Khasminskii-type condition parameter $p$, super-linear growth constant $ \g$, repectively.
%
Letting $ p \to \infty$, then \eqref{rate22} gives $$  R \to \frac{1}{1 + \g}.$$ This is almost the optimal $\mathcal L ^r$-convergence rate of the truncated EM method in the case of jumps. Only in the case of $\g =0$, i.e. the drift and diffusion coefficients grows linearly, this convergence rate is close to 1.
It should be mentioned that this is significantly different from the result on SDEs without jumps. We already known that for any $ r \ge 2$ (see \cite{Guo2018note})
$$ \E | x(T) - x_{\DD}(T)|^r \le C \DD ^ { r ( 1- 2 \veps)/2} ,\qu \forall \veps \in (0, 1/4],$$
which means that the $\mathcal L^r$-convergence rate is close to $r/2$ when there is no jumps in SDE \eqref{eq0}. In fact, this difference is caused by the following reason: all moments of the Poisson increments $\DD N_k = N( (k+1)\DD) - N( k \DD)$ have the same order $ \DD$ (see \eqref{memont of Nk}), while the Brownian increments $\DD B _k = B( (k+1)\DD) - B( k \DD) $ have different orders, namely $ \E |\DD B _k|^ {2m} = o (\DD ^m)$ and $ \E |\DD B_k|^ {2m+1} =0$. These properties force the control function $ \ph (\DD) := \DD ^ {- \veps}$ to be not greater than $ \DD ^ {-1/p} \we \DD ^ {-1/4}$, when we are tring to bound the moments of the truncated EM solution in Lemma \ref{moment of xdt}. This eventually leads to the differences in the convergence rates between SDEs with and without jumps.
\end{remark}

\subsection{Convergence and convergence rate of the truncated EM method in $ \mathcal L^r( 0 <  r  < 2)$}\label{lessthan2}
In this subsection, we will discuss the convergence in $\mathcal L^r( 0 <  r  < 2)$ under the assumption that the drift, diffusion and jump terms behave like a polynomial. For this purpose, we first impose the following assumptions.
\begin{assumption} \label{2Local_condition}
There exists positive constant $ K_n $ such that
\begin{align} \label{Local_eq}
 \left| {f(x) - f(y)} \right| \vee  {\left | {g(x) - g( y)} \right|} \vee  {\left | {h(x) - h( y)} \right|}  \le K_n{| {x - y} |},  \qu \forall x,y \in \R ^d, \; |x| \ve |y| \le n .
\end{align}
\end{assumption}

\begin{assumption} \label{2monotone_condition}
There exist constants $ \bar K>0$ such that
\begin{align}\label{monto_eq}
 2 x^ T f(x)   +  | g(x)|  ^2  + \la ( 2 x^ T h(x) + |h(x)|^2) \le \bar K(1 + |x|^2 ),  \quad  \forall x \in \R ^d.
\end{align}
\end{assumption}
We also give a known result as a lemma.
\begin{lemma} \label{2lemma4.77}
Under Assumption \ref{2Local_condition} and \ref{2monotone_condition} the SDE \eqref{eq0} has a unique global solution $x(t)$, moreover,
\begin{align}\label{}
\sup _ {0 \le t \le T} \E |x(t)| ^ 2 < \infty , \quad \forall T > 0.
\end{align}
\end{lemma}
In this subsection, all the three coefficients of the SDEs are allowed to grow super-linearly. Hence, we have to truncate the three terms. Similarly, we first choose a strictly increasing function $\mu : \R ^+  \rightarrow \R^ +$ such that $\mu (n) \rightarrow \infty $, as $ n \rightarrow \infty$,
and
\begin{align}\label{2HD}
\sup _{|x| \le n}  |f(x)| \vee |g(x)| \ve |h(x)| \le \mu (n), \quad \forall n \ge 1.
\end{align}
Denoted by $\mu ^{-1}$ is the inverse function of $\mu$. We also choose a constant $\Dst \in (0, 1 )$ and a strictly decreasing function $\ph: (0, \Dst) \rightarrow (0, \infty)$ such that
 \begin{align}\label{2delta_h_relations}
 \ph( \Delta ^{\ast})  \ge \mu (1), \quad \lim _ {\Delta \to 0} \ph(\Delta) = \infty \quad  \textrm{  and} \quad  \ph(\DD) \DD ^ {1/4} \le 1  .
 \end{align}
For a given step size $\Delta \in (0, \Delta ^ {\ast})$, the truncated functions are defined as below
\begin{align*}
f_{\Delta}(x) = f( \pi _ {\Delta}(x)) , \quad
 g_{\Delta}(x) = g( \pi _ {\Delta}(x))\quad \textrm{and} \quad h_{\Delta}(x) = h( \pi _ {\Delta}(x)), \quad \forall x \in \R^d,
\end{align*}
where $\pi _{\Delta }$  is defined as the same as before. The following lemma also shows that the truncated functions preserve the Khaminskii-type condition. The proof is given in the Appendix.
\begin{lemma} \label{2lem 24}
Let Assumption \ref{2monotone_condition} hold. Then, for all $\Delta \in (0, \Delta ^ \ast ]$, we have
\begin{align}\label{ddd}
2 x^{T} f_ {\Delta} (x) +  |g_{\Delta}(x)|^2   +  \la ( 2 x^T h_{\DD}(x) + |h_{\DD} (x)|^2 )\le 2 \bar K (1+ |x|^2), \quad \forall x \in \R^d.
\end{align}
\end{lemma}
We now give the discrete-time truncated EM scheme
\begin{align}\label{2EM scheme}
\XD (t_{k+1}^- ) = \XD (t_k^-) + \fD(\XD(t_k^-)) \DD + \gD (\XD (t_k^-)) \DD B_k + \hD (\XD (t_k^-)) \DD N_k, \quad 0 \le k \le M-1,
\end{align}
where $X_{\Delta}(t_k^-)$, $\DD B_k $ and  $\DD N_k $ is defined the same as before. The continuous-time form is defined by
\begin{align}\label{}
\xdt  = x_0 + \int_0^t \fdxs ds + \int_0^t \gdxs dB(s) + \int_0^t \hdxs dN(s).
\end{align}
where $ \bxdt$ is defined in \eqref{step_process}. In order to state our main result, we first give some useful lemmas.
\begin{lemma} \label{2lem17}
For any $\DD \in (0, \Delta ^ {\ast})$ and $t >0$, we have
\begin{align}\label{2eq114}
\E | \xdt  - \bxdt|^{\hat p} & \le C_{\hat p}   (\ph(\DD))^{\hat p} \DD, \quad  {\hat p} \ge 2, \\
\E | \xdt  - \bxdt|^{\hat p} & \le C_{\hat p}   (\ph(\DD))^{\hat p} \DD ^ {\hat p /2}, \quad   0 < {\hat p} < 2.
                                                \end{align}
\end{lemma}
 Consequently,
\begin{align}\label{2eq115}
\lim _ {\DD \to 0} \E |\xdt - \bxdt|^{\hat p} =0 , \quad \forall t \ge 0.
\end{align}
\pr  Fix any $\DD \in (0, \DD ^ {\ast}]$, $t \ge 0$ and $\hat p \ge 2$. There is an integer $k \ge 0$ such that $\tk \le t < \tkk$. By \As, we have
\begin{align}\label{2eq116}
& \E | \xdt - \bxdt|^{\hat p}  \\ \nonumber
 & \le C_{\hat p} \left ( \E \left | \int_{\tk}^t \fdxs ds \right |^{\hat p} + \E \left | \int_{\tk}^t \gdxs dB(s) \right|^{\hat p}
    + \E \left | \int_{\tk}^t \hdxs dN(s) \right |^{\hat p}      \right )\\ \nonumber
 & \le C_{\hat p} \Big ( \DD ^ {\hat p-1} \E  \int_{\tk}^t\left | \fdxs  \right |^{\hat p} ds + \DD ^ { (\hat p -2)/2
 }\E  \int_{\tk}^t \left | \gdxs  \right|^{\hat p} ds +
                      \E \left | \int_{\tk}^t \hdxs d N (s) \right |^{\hat p}   \Big ) \\ \nonumber
 &  \le C_{\hat p} \Big (\DD ^ {p/2} (\ph(\DD))^{\hat p}   +
                      \E \left | \int_{\tk}^t \hdxs d N (s) \right |^{\hat p}   \Big ),
\end{align}
where $C_{{\hat p}}$ is a generic constant.
The property of Poisson increments implies
\begin{align*}
\E \left | \int_{\tk}^t \hdxs d N (s) \right |^{\hat p}  & \le (\ph(\DD))^{\hat p} \E |\DD N_k|^{\hat p} \no
                              & \le c_0 (\ph(\DD))^{\hat p} \DD.
        \end{align*}
Inserting this into \eqref{2eq116} and recalling $\hat p \ge 2$ gives
\begin{align*}
\E | \xdt - \bxdt |^{\hat p} \le C_{\hat p} (\ph(\DD))^{\hat p} \DD.
\end{align*}
Noting from \eqref{delta_h_relations} that $(\ph(\DD))^{\hat p} \Delta  = (\ph(\DD))^{\hat p} \Delta ^ {1/2} \Delta ^ {1/2}  \le \DD ^ {1/2}$, we obtain \eqref{2eq115} form \eqref{2eq114}. \\
For $ 0 < \hat p < 2$, we have
\begin{align*}\label{}
& \E | \xdt - \bxdt |^{\hat p} \le \Big ( \E | \xdt - \bxdt |^2  \Big ) ^ {\hat p /2} \no
& \le \Big(  C_{\hat p} (\ph(\DD))^{2} \DD \Big ) ^ {\hat p /2}  =  C_{\hat p}(\ph(\DD))^{\hat p } \DD ^ {\hat p /2}.
\end{align*}
Thus, the proof is complete. $\Box$

The following lemma give the moment bound of the truncated EM solutions.
\begin{lemma} \label{2moment of xdt}
Let \As \ref{2Local_condition} and \ref{2monotone_condition} hold. Then
\begin{align}\label{2eq moment lof xdt }
\sup _{ 0 \le \DD \le \Dst} \sup _ {0\le t \le T}  \E | \xdt|^2 \le C, \qu \forall T > 0.                                                                     \end{align}
\end{lemma}
\pr Fix any $ \DD \in (0,\Dst )$ and $T >0 $. By the \Ito and \As \ref{2monotone_condition}, we have
\begin{align}\label{2eq119}
  \E | \xdt|^2  & \le \E |x_0|^2  + \E \int_0^t \Big ( 2  x^T _ {\DD} (s^-) \fdxs  +  |\gdxs|^2  \Big )ds \no
   & \qu + \lambda \E  \int_0^t \Big (  2 x_{\DD}(s^-)^T \hdxs + | \hdxs |^2  \Big )ds \no
& \le \E |x_0|^2  + \E \int_0^t \Big ( 2  \bar x _ {\DD}  ^ T(s^-) \fdxs  +  |\gdxs|^2  \Big )ds \no
   & \qu + \lambda \E  \int_0^t \Big (  2 \bar x _ {\DD}  ^ T(s^-)\hdxs + | \hdxs |^2  \Big )ds + \bar J_1 \no
   & \le  \E |x_0|^2 + 2 \bar K \int_0^t \Big ( 1 + \E | \bxds |^2 \Big )ds  + \bar J_1,
\end{align}
where
\begin{align*}
\bar J_1 =  \E \int_0^t \Big ( 2 ( x_{\DD}(s^-)  -  \bar x _ {\DD}(s^-))^T \fdxs  +  2 \la ( x_{\DD}(s^-)  -  \bar x _ {\DD}(s^-))^T \hdxs \Big )ds.
\end{align*}
By Lemma \ref{2lem17}, \eqref{2HD} and \ref{2delta_h_relations}, we have
\begin{align*}
\bar J_1  & \le   2(\la + 1) \ph(\DD)  \int_0^t  \E |   x_{\DD}(s^-)  -  \bar x _ {\DD}(s^-)| ds \no
& \le    2(\la + 1) T C ( \ph(\DD)) ^2 \DD ^ {1/2} \le    C .
\end{align*}
Inserting this into \eqref{2eq119} and using Lemma \ref{2lem17} gives
\begin{align*}
\E | \xdt|^2  \le C  + 2 \bar K \int_0^t  \E  | \bxds |^2  ds .
\end{align*}
Hence, we have
\begin{align*}
\sup _ {0 \le u \le t }\E | x_{\DD}(u) |^2  \le C  + 2 \bar K \int_0^t \sup _ {0 \le u \le s }\E | x_{\DD}(u) |^2   ds .
\end{align*}
The Gronwall inequality yields
$$  \sup _ {0 \le u \le T} \E  |  x _ {\DD} (u)|^p \le C.$$
Thus, we complete the proof. $\Box$

As the proof is in a similar way as Lemma \ref{lemma34} and \ref{lemma3.2} were proved, we also have the following Lemma.
\begin{lemma} \label{2lemma33}
Let \As \ref{2Local_condition} and \ref{2monotone_condition} hold. For any real number $n > |x_0| $,
              then
              \begin{align}\label{eq35}
\PP (\tr \le T)  \le  \frac{C}{n^2} \qu \textrm{and} \qu \PP (\rho _{\DD , n } \le T)  \le  \frac{C}{n^2},
              \end{align}
where $ \tr$ and $ \rho _{\DD , n } $  is the same as before.

\end{lemma}
Now, let us show the convergence of the truncated EM methods for the SDEs.
 \begin{theorem}\label{2main_result34}
Let \As \ref{2Local_condition} and \ref{2monotone_condition} hold. Then, for any $r \in (0,2 )$
\begin{align}\label{re1}
\lim _ { \DD \to 0 }\E |x(T) - x_{\DD}(T)| ^r  =0
\end{align}
and
\begin{align}\label{re2}
\lim _ { \DD \to 0 } \E |x(T) - \bar x_{\DD}(T)| ^r =0 .
\end{align}
\end{theorem}
\pr Let $\tr$, $\rdr$, $\tdr$ and $ e _ {\DD}(t) $ be the same as before. Applying the Young inequality, we have that for any $\delta >0 $,
\begin{align}\label{2eq330}
\E |e_ {\DD } (T)|^r  & = \E \Big ( | e_{\DD}(T)|^r \II _ {\{ \tdr > T \}  }\Big ) +  \E \Big ( | e_{\DD}(T)|^r \II _ {\{ \tdr \le T  \}  }\Big ) \no
  & \le \E \Big ( | e_{\DD} (T \we \tdr )|^r \Big ) +  \frac{r \delta}{ 2} \E |e_ {\DD } (T)|^2 + \frac{2-r}{2 \delta ^ {r/ (2-r)}}\PP ( \tdr \le T).
\end{align}
By Lemma \ref{2lemma4.77} and \ref{2moment of xdt}, we have
\begin{align}\label{2eq331}
\E |e_ {\DD } (T)|^2 \le 2 \E |x (T)|^p + 2 \E |x_ {\DD } (T)|^p \le C.
\end{align}
Using Lemma \ref{2lemma33}, we obtain
\begin{align}\label{2eq332}
\PP ( \tdr \le T) \le \PP ( \tau_n \le T) + \PP (\rdr \le T ) \le \frac{C}{n^2}.
\end{align}
Inserting \eqref{2eq331} and \eqref{2eq332} into \eqref{2eq330}, we get
\begin{align*}
\E |e_ {\DD } (T)|^r  \le \E | e_{\DD} (T \we \tdr )|^r + \frac{Cr \delta}{2} + \frac{C(2-r)}{2 n^2 \delta ^ {r/(2-r)}}.
\end{align*}
Now, let $\veps >0$ be arbitrary. We can choose $\delta$ sufficiently small such that $$ \frac{Cr \delta}{2} \le \frac{\veps}{3} $$ and then choose $n$ sufficiently large such than  $$ \frac{C(2-r)}{2 n^2 \delta ^ {r/(2-r)}} \le \frac{\veps}{3}.$$
We may assume that $\Dst $ is sufficiently small for $ \mu ^{-1} (\ph(\Dst)) \ge n$. In the same way as Theorem 3.5 in \cite{Mao2015} was proved, we can show that for all $\DD \in (0, \Dst]$
$$ \E |e_ {\DD } (T)|^2 \le C \DD ,$$
which implies  $$ \E \Big ( | e_{\DD} (T \we \tdr )|^r \Big ) \le  \frac{\veps}{3} .$$
Hence, we obtain the required assertion \eqref{re1}. Combining this with Lemma \eqref{2lem17} gives \eqref{re2}.
Thus, the proof is complete. $\Box$


For the purpose of getting the convergence rates a time T, we need some additional conditions.
\begin{assumption} \label{2onesided_conditions}
 There exists a constant $\bar L_1 >0$  such that
\begin{align} \label{2onesided_eq}
& 2 (x-y)^T (f(x) -f(y) )  + |g(x) - g(y)|^2  \no
& + 2 \la  (x-y)^T (h(x) -h(y) ) +  \la |h(x) - h(y)|^2\le \bar L_1 |x -y | ^2,
\end{align}
for any $ x,y \in \R ^d$.
\end{assumption}
\begin{assumption} \label{2Polynomial_condition}
There exist constant $\bar L_2 >0$ and $  \bar \g  \ge 0$ such that
\begin{align} \label{2Polynomial_eq}
|f(x) - f(y)| \ve |h(x) - h(y)| \le \bar L_2(1+ |x|^{\bar \g} + |y|^{\bar \g}) |x-y|, \quad \forall  x , y \in \R ^d.
\end{align}
\end{assumption}
Obviously, this condition implies
\begin{align}\label{L4444}
|f(x)| \ve |h(x)| \le \bar L_3 |x| ^ {1+ \bar \g},
\end{align}
where $ \bar L_3  = 2 \bar L_1 + |f(0)| + |g(0)|$.
\begin{lemma} \label{2lemma3.2}
Let \As \ref{2Local_condition}, \ref{2monotone_condition}, \ref{2onesided_conditions} and \ref{2Polynomial_condition} hold.
Let $n > |x_0| $ be a real number and let $\DD \in (0, \Dst ) $ be sufficiently small such $\muh \ge n $. Let $\tr $ and $\rdr$ be the same as before. Let
\begin{align*}
\tdr = \tr \we \rdr \qu \textrm{and } \qu   \eds = x(t^-) - x_{\DD} (t^-) \qu \forall t> 0.
\end{align*}
Then  \begin{align*}
\E | e _ {\DD} (T \we \tdr ) |^2 \le C  (\ph(\DD))^ {2} \DD .
       \end{align*}
\end{lemma}
\pr We write $\tdr = \theta$ for simplicity. By the \Ito and \As \ref{2onesided_conditions}, we get that for $ 0 \le t \le T$,
\begin{align}\label{2eq43}
& \E |e _ {\DD} (t \we \theta)|^2   \no
& \le \E \int_0^{t \we \theta}   \Big (  2e^T _ {\DD}(s^-) (f(x(s^-) )- \fdxs) +  |g(x(s^-))  -\gdxs | ^2 \Big ) ds\no
&  \qu + \la \E \int_0^{t \we \theta}    \Big ( |\eds + (h(x(s^-))- \hdxs) |^2  - | \eds |^2 \Big )ds \no
& \le \E \int_0^{t \we \theta}   \Big (  2(x(s^-)- \bar x _ {\DD}(s^-)  )^T  (f(x(s^-) )- \fdxs) +  |g(x(s^-))  -\gdxs | ^2 \Big ) ds  + \bar J_2\no
& \qu +  \E \int_0^{t \we \theta}   \Big (  2 \la (x(s^-)- \bar x _ {\DD}(s^-)  )^T  (h(x(s^-) )- \hdxs) +  \la |h(x(s^-))  -\hdxs | ^2 \Big ) ds  + \bar J_3 \no
& \le \bar L_1   \int_0^{t } \E |  x(s \we \theta )- \bar x_{\DD}(s\we \theta  )| ^2 ds  + \bar J_2 + \bar J_3,
\end{align}
where \begin{align*}
\bar J_2 = 2 \E \int_0^{t \we \theta}     (x_{\DD }(s^-)- \bar x _ {\DD}(s^-)  )^T  (f(x(s^-) )- \fdxs)  ds , \no
\bar J_3 = 2 \la  \E \int_0^{t \we \theta}     (x_{\DD }(s^-)- \bar x _ {\DD}(s^-)  )^T  (h(x(s^-) )- \hdxs)  ds .
      \end{align*}
By the condition $\muh \ge n $ and the definition  of the truncated functions $\fD$ and $\gD$, we have that
\begin{align*}
\fdxs = f(\bxds) \qu \textrm{and} \qu \gdxs = g(\bxds), \; \textrm{for} \; 0 \le s \le t \we \theta.
\end{align*}
%
Hence, by \As \ref{2Polynomial_condition} and the \Holder as well as Lemma  \ref{2lem17} and  \ref{2moment of xdt},  we get that
\begin{align}\label{2eq1}
\bar J_2 & \le  2 \E \int_0^{t \we \theta}     | x_{\DD }(s^-)- \bar x _ {\DD}(s^-)  ||  f(x(s^-))- f(\bxds)|ds  \no
& \le  2\bar L_2 \E \int_0^{t \we \theta}| x_{\DD }(s^-)- \bar x _ {\DD}(s^-)  ||1+ |x(s^-)|^{\bar \g}+ |\bxds|^{\bar \g}||  x(s)- \bxds|ds \no
& \le  \bar L_2 \int_0^{t \we \theta} \E |  x(s^-)- \bxds| ^2 ds + C \int_0^{t \we \theta} \E (1+ |x(s^-)|^{2 \bar \g}+ |\bxds|^{2 \bar \g})|  x_{\DD}(s^-)- \bxds|^2 ds \no
& \le  \bar L_2 \int_0^{t } \E |  x(s \we \theta )- \bar x_{\DD}(s\we \theta  )| ^2 ds \no
 & + C \int_0^{T } \B (1+ \E |x(s^-)|^{2}+ \E |\bxds|^{2}\B )^ {\bar \g }  \B( \E  | x_{\DD}(s^-)- \bxds|^{2/(1 - \bar
 \g )}  \B )^{1- \bar \g}ds \no
 & \le  \bar L_2 \int_0^{t } \E |  x(s \we \theta )- \bar x_{\DD}(s\we \theta  )| ^2 ds  + C (\ph(\DD))^2 \DD .
\end{align}
Similarly, we have
\begin{align}\label{2eq2}
\bar J_3
 & \le  \la  \bar L_2 \int_0^{t } \E |  x(s \we \theta )- \bar x_{\DD}(s\we \theta  )| ^2 ds  + C  (\ph(\DD))^2 \DD.
\end{align}
Inserting \eqref{2eq1}, \eqref{2eq2} into \eqref{2eq43} and combining Lemma \ref{2lem17}, we have
\begin{align*}
 \E |e _ {\DD} (t \we \theta)|^2  \le C \intt \E | e_{\DD} (s \we \theta )|^2 ds +C  (\ph (\DD))^ {2}\DD.
\end{align*}
The Gronwall inequality complete the proof. $\Box$

 \begin{theorem}\label{2main_result4}
Let \As \ref{2Local_condition}, \ref{2monotone_condition}, \ref{2onesided_conditions} and \ref{2Polynomial_condition} hold. Let $r \in (0, 2)$.
If
\begin{align}\label{2cond3}
\ph(\DD) \ge  \mu \Big ( \bar L_3 ^ {-(1+ \bar \g)}  ( (\ph(\DD))^ {r}\DD ^ {r /2 }) ^ {-1 /{(2-r)}}\Big )
   \end{align}
holds for all sufficiently small $\DD \in (0, \Dst)$, then for every such small $\DD$,
\begin{align}\label{}
\E |x(T) - x_{\DD}(T)| ^r  \le C (\ph(\DD))^ {r}\DD^{r/2}
\end{align}
and
\begin{align}\label{2eq58}
\E |x(T) - \bar x_{\DD}(T)| ^r \le C  (\ph(\DD))^ {r}\DD^{r/2} .
\end{align}
\end{theorem}
\pr Let $\tr$, $\rdr$, $\tdr$ and $ e _ {\DD}(t) $ be the same as before. By \eqref{2eq330} - \eqref{2eq332}, inequality
\begin{align*}
\E |e_ {\DD } (T)|^r  \le \E | e_{\DD} (T \we \tdr )|^r + \frac{Cr \delta}{2} + \frac{C(2-r)}{2 n^2 \delta ^ {r/(2-r)}}
\end{align*}
holds for any $\DD \in (0 , \Dst)$, $n > |x_0|$ and $\delta>0$.  We can therefore choose $\delta = (\ph(\DD))^ {r} \DD ^ {r/2} $ and $n = \bar L_3 ^ {- (1+ \bar \g)}( (\ph(\DD))^ {r} \DD ^ {r/2})^ {- 1/{(2-r)}}$ to get
\begin{align*}
\E |e_ {\DD } (T)|^r  \le \E | e_{\DD} (T \we \tdr )|^r + C  (\ph(\DD))^ {r}\DD ^ {r/2}.
\end{align*}
By condition \eqref{2cond3}, we have
\begin{align*}
\mu ^{-1} (\ph(\DD)) \ge  \bar L_4 ^ {- {1 + \bar \g }} ((\ph(\DD))^ {r} \DD ^ {r/2} )^ {-1 /{(2-r)}} =n.
\end{align*}
Using Lemma \ref{2lemma3.2}, we have
\begin{align*}
\E |e_ {\DD } (T)|^r  \le  (\E |e_ {\DD } (T)|^2 )^{r/2}  \le    C (  (\ph(\DD))^ {2} \DD )^{r/2} =  C  (\ph (\DD))^ {r}\DD ^ {r/2}.
\end{align*}
Combining this with Lemma \eqref{2lem17} gives \eqref{2eq58}. Thus, the proof is complete. $\Box$
\begin{corollary}\label{cor2}
Let \As \ref{2Local_condition}, \ref{2monotone_condition}, \ref{2onesided_conditions} and \ref{2Polynomial_condition} hold. Define
\begin{align}\label{re444}
\mu (n ) = \bar L_3 n ^ {1+ \bar \g}, \qu n \ge0.
\end{align}
Let $0 < r \le 2/(2 + \bar \g)$ and
\begin{align}\label{re333}
\ph (\DD ) = \DD ^ {- \veps}, \qu \veps \in \B [\frac{r (1 + \bar \g)}{ 4+ 2 r  \bar \g}, \frac{1}{4}\B ].
\end{align}
Assume that  $\bar \DD \in (0,1)$ and
\begin{align}\label{2cond4}
\ph(\bar \DD) \ge  \mu \Big ( \bar L_3 ^ {-(1+ \bar \g)}  ({\bar \DD } ^ {r /2 } (\ph({\bar \DD}))^ {r}) ^ {-1 /{(2-r)}}\Big )
   \end{align}
hold. Then, for every $\DD \in (0, \bar \DD)$ ,
\begin{align}\label{re111}
\E |x(T) - x_{\DD}(T)| ^r  \le C \DD^{r/2 - r \veps}
\end{align}
and
\begin{align}\label{re222}
\E |x(T) - \bar x_{\DD}(T)| ^r \le C \DD^{r/2 - r \veps }  .
\end{align}

\end{corollary}
\pr Applying Theorem \ref{2main_result4} along with \eqref{re444} and \eqref{re333} gives the required assertion
\eqref{re111} and \eqref{re222}. $\Box$

\begin{remark}
Substituting \eqref{re444} and \eqref{re333} into \eqref{2cond4} gives
 \begin{align*}
 \bar \DD ^ {- \veps} \ge \bar \DD ^ {   - r (1 - 2 \veps) (1 + \bar \g)/ (2 (2 - r)) }, \qu \textrm{namely} \qu \veps  \ge \frac{r (1 + \bar \g)}{ 4+ 2 r  \bar \g}.
 \end{align*}
But, condition \eqref{re333} means
\begin{align*}
\frac{r (1 + \bar \g)}{ 4+ 2 r  \bar \g} \le \frac{1}{4}, \qu \textrm{namely} \qu r \le \frac{2}{ 2 + \bar \g} \le 1.
\end{align*}
Hence, we have to force $r$ to be not greater than $2/ (2 + \bar \g)$ in the corollary \ref{cor2}.
\end{remark}
\begin{remark}
Fixing $\bar \g \ge 0$, by \eqref{re333} and \eqref{re111}, we can conclude that convergence rate is a increasing function with respect to $\veps$. Hence, substituting $$ \veps = \frac{ r (1 + \bar \g)}{ 4 + 2 r \bar \g}$$ into $r/2( 1 -2 \veps) $ obtains the optimal $ \mathcal L^r$-convergence rate, that is
\begin{align}\label{rate222}
 R := \frac{ r (2-r)}{ 2 ( 2 + r \bar \g)} ,
\qu \textrm{ for} \qu 0 < r \le \frac{2}{ 2 + \bar \g },
\end{align}
which means convergence rate $R$ increases as $r$ increases. In other words, the higher moment has a better convergence rate for SDEs with jumps when $ 0< r \le 2/(2 + \bar \g)$. If we take $$r = \frac{2}{2+ \bar \g},$$ then \eqref{rate222} becomes
$$ R = \frac{1}{4 +  2 \bar \g }, $$
this is the maximum of optimal $\mathcal L^r$-convergence rate.
In particular, if $\bar \g =0$, i.e. the drift and diffusion coefficients grows linearly, then convergence rate is equal to $1/4$  by choosing $r=1$.
%
\end{remark}

\section{Asymptotic behaviours}
\subsection{Stability}\label{secstab}
In this subsection, we will show that the partially truncated EM method can preserve the mean square exponential stability of the underlying SDEs \eqref{eq0}.
For the purpose of stability, we also assume that
\begin{align}\label{eq4.1}
f(0) = g(0) = h(0)  = 0,
\end{align}
which means
\begin{align}\label{hk0}
|F_1(x)| \ve |G_1(x) | \ve | h(x)|\le K_1 |x| ,\qu \forall x \in \R^d.
\end{align}
We first impose the following assumption.
\begin{assumption} \label{stability_condition}
Assume that there exist positive constant $\a_1 \ge \a_2 + \la K_1(2 + K_1)$  and $\theta \in (0, \infty)$ such that
\begin{align*}
2x^T F_1(x) + (1+\theta) |G_1(x)|^2 \le - \a_1 |x|^2
\end{align*}
and
\begin{align*}
2x^T F(x) + (1+\theta ^ {-1}) |G(x)|^2 \le  \a_2 |x|^2
\end{align*}
for all $x \in \R ^d$.
\end{assumption}
If there is no super-linear term $G(x)$, we set $\theta = 0$ and $ \theta ^{-1} |G(x)|^2 =0$. Similarly, when linear term $G_1(x)$ is absent, we set $\theta  =  \infty $ and $ \theta  |G_1(x)|^2 =0$. Moreover, this assumption means \begin{align}\label{}
2x^T f(x) + |g(x)|^2 + \la ( x^T h(x) + |h(x)|^2) \le -(\a_1 - \a_2 - \la K_1 (2 + K_1))|x|^2, \qu x \in \R^d.
                      \end{align}
It is therefore known that the SDEs \eqref{eq0} is exponentially stable in the mean square sense. We state the following Lemma.
\begin{lemma}\label{stability_theorem}
Let \As \ref{Local_condition} - \ref{monotone_condition} and \ref{stability_condition} hold. Then for any initial value $x_0 \in \R^d$, the solution of the SDEs \eqref{eq0} satisfies
 \begin{align*}\label{}
 \E |x(t)|^2 \le |x_0|^2 e^{-(\a_1 - \a_2 - \la K_1 (2 + K_1))t^-}  ,\qu \forall   t \ge 0.
 \end{align*}
\end{lemma}

The following theorem shows that the truncated EM method preserves this mean square exponential stability perfectly. We will employ the technique due to Guo et al. \cite{Guo2017partial} to prove our results.
\begin{theorem}\label{stability_theorem2}
Let \As \ref{Local_condition}, \ref{onesided_conditions}, \ref{monotone_condition} and \ref{stability_condition} hold. Then for any $ \epsilon \in (0,\a_1 - \a_2 - \la K_1 (2 + K_1)) $, there exists a $\hat \DD
\in (0, \Dst) $ such that for all $ \DD \in (0, \hat \DD)$ and any initial value $x_0 \in \R ^d$, the solution of the truncated EM method \eqref{EM scheme} satisfies
\begin{align}\label{eq401}
\E |\XD (\tk)|^2 \le |x_0|^2 e ^{-(\a_1 - \a_2 - \la K_1 (2 + K_1) - \epsilon) \tk  }, \qu \forall  k \ge 0.
\end{align}
\end{theorem}
\pr Fix $\DD \in (0, \Dst)$. In the same way as Theorem 4.3 in \cite{Guo2017partial} was proved,  we have
\begin{align}\label{eq402}
2x^T f_{\DD} (x) + |g _ {\DD}(x)|^2 \le - (\a_1 - \a_2  ) |x|^2, \qu \forall x \in \R^d.
\end{align}
From \eqref{EM scheme}, we have
\begin{align}\label{eq406}
\E | \XD ( \tkk) |^2 &  = \E \Big ( |\xdtk|^2 + |\fD (\xdtk)| ^2 \DD ^2   + |\gD (\xdtk) \DD B_k|^2   \no
& \qu  + 2 \xdtk ^T  \fD (\xdtk)  \DD  +  |h (\xdtk) \DD N_k|^2   \no
& \qu +  2 \DD \fD ^T (\xdtk) h (\xdtk)\DD N_k   +  2 \xdtk ^T h (\xdtk)\DD N_k    \Big )
\end{align}
for $ 0 \le k  \le M-1.$
The property of Brownian increments implies
\begin{align*}
\E |\gD (\xdtk) \DD B_k|^2  = \DD \E |\gD (\xdtk) |^2 .
\end{align*}
But, the Poisson increments satisfy $\E \DD N_k = \la \DD $ and $\E (\DD N_k)^2 = \la \DD ( 1+ \la \DD)$.
Hence, using the independence of the increments and \eqref{hk0}, we find that
\begin{align}\label{eq403}
2 \E | \xdtk  h (\xdtk)\DD N_k | \le 2  K_1 \E | \xdtk|^2 \E | \DD N_k| =  2  K_1 \la \DD  \E | \xdtk|^2,
\end{align}
\begin{align}\label{eq404}
 \E |h (\xdtk) \DD N_k|^2  & \le K_1^2 \E |\xdtk|^2 \E | \DD N_k|^2 \no
  & \le K_1 ^2 \la \DD (1+ \la \DD) \E |\xdtk|^2 \no
  & = K_1 ^2 \la \DD  \E |\xdtk|^2 + K_1 ^2 \la ^2 \DD ^2 \E |\xdtk|^2
\end{align}
and
\begin{align}\label{eq405}
2 \E |\DD \fD  (\xdtk) h (\xdtk)\DD N_k |  & \le 2 K_1 \DD \E ( |\xdtk \fD (\xdtk) | ) \E |\DD N_k| \no
 & \le  K_1 \la \DD ^2 (  \E | \xdtk|^2 + \E | \fD (\xdtk)| ^2 )  .
\end{align}
Substituting \eqref{eq403}-\eqref{eq405} into \eqref{eq406} gives
\begin{align}\label{eq408}
\E | \XD ( \tkk) |^2 &  \le  \E \Big ( |\xdtk|^2 +  2 \xdtk ^T  \fD (\xdtk)  \DD  + |\gD (\xdtk) |^2   \DD  \Big )  \no
& \qu + \la K_1 (2 + K_1) \DD  \E | \xdtk |^2  +  (1+ K_1 \la )  \DD ^2  \E |\fD (\xdtk)| ^2  \no
& \qu + (K_1^2 \la ^2 + K_1 \la) \DD^2 \E | \xdtk |^2.
\end{align}
By \eqref{eq402}, we have
\begin{align}\label{eq409}
\E | \XD ( \tkk) |^2 &  \le  (1- (\a_1- \a_2 - \la K_1 (2 + K_1)  )\DD) \E |\xdtk|^2   \no
& \qu  +  (1+ K_1 \la )  \DD ^2  \E |\fD (\xdtk)| ^2   + (K_1^2 \la ^2 + K_1 \la) \DD^2 \E | \xdtk |^2.
\end{align}
By \eqref{Polynomial_eq} and \eqref{eq4.1}, we have
\begin{align*}
|F_{\DD} (x)| ^2 \le 4 L_1 |x| ^2, \; \textrm{if} \; |x| \le 1,
\end{align*}
and
\begin{align*}
|F_{\DD} (x)| ^2 \le  (\ph( \DD ))  ^2 \le (\ph( \DD ))  ^2 |x|^2, \; \textrm{if} \; |x| > 1.
\end{align*}
Hence, we have
\begin{align*}
\DD |f_{\DD} (x)| ^2 & \le 2 (K_1^2+ 4 L_1+ (\ph( \DD ))  ^2 ) \DD |x|^2 \no
 &  \le  2 \Big (  (K_1^2+ 4 L_1) \DD + \DD ^ { 1/2  \we (\bar p -2) / \bar p   } \Big )   |x|^2
\end{align*}
for all $x \in \R^d$, where \eqref{delta_h_relations} have been used.
For any $ \epsilon \in (0, \a_1 - \a_2 - \la K_1 (2 + K_1))$, there is a $ \hat \DD \in (0, \Dst)$ sufficiently small such that for all $\DD \in (0, \hat \DD)$, $ ( \a_1 - \a_2 - \la K_1 (2 + K_1) ) \DD \le 1 $ and
\begin{align}\label{eq4.2}
\begin{cases}
2(1 + K_1 \la) \Big (  (K_1^2+ 4 L_1) \DD + \DD ^ { 1/2  \we (\bar p -2) / \bar p  } \Big )   \le 0.5 \epsilon \\
(K_1^2 \la ^2 + K_1 \la ) \DD  \le 0.5 \epsilon.
 \end{cases}
\end{align}
For each such $\DD$, we have
\begin{align*}
(1+ K_1 \la )  \DD ^2  \E |\fD (\xdtk)| ^2   + (K_1^2 \la ^2 + K_1 \la) \DD^2 \E | \xdtk |^2 \le \epsilon \DD \E | \xdtk |^2.
\end{align*}
Inserting this into \eqref{eq409}, we yield
\begin{align}\label{eq411}
\E |\XD ( \tkk) |^2 &  \le  (1- (\a_1 - \a_2 - \la K_1 (2 + K_1)  - \epsilon )) \DD )  \E |\xdtk|^2   \no
                    & \le \cdots \no
                 & \le |x_0|^2(1- ( \a_1 - \a_2 - \la K_1 (2 + K_1) - \epsilon ) ) \DD )^{k+1}.
\end{align}
Elementary inequality  $1 + x \le e ^x ,\; \forall x \in \R$ gives
\begin{align}\label{eq412}
\E | \XD ( \tkk) |^2   \le  |x_0|^2  e ^ {- (\a_1 - \a_2 - \la K_1 (2 + K_1)  - \epsilon ) t_{k+1}^- }.
\end{align}
Thus, the proof is complete.  $\Box$

\subsection{Asymptotic boundedness} \label{secbound}
In this  subsection, we will show that the truncated EM method maintains the asymptotic boundedness of the underlying of SDEs \eqref{eq0}. The additional assumption is the following one.

\begin{assumption} \label{boundedness}
Assume that there exist constant $\bar \alpha_1 ,\bar \a _2 , \bar \b_2 > 0$ and $\bar \b_1> \bar \b_2 + \max(\la (4K_1^2 +1), 2 \la K_1 (2 + K_1))$ such that
\begin{align*}
2x^T F_1(x) + (1+\theta) |G_1(x)|^2 \le \bar \a_1 - \bar \b_1 |x|^2
\end{align*}
and
\begin{align*}
2x^T F(x) + (1+\theta ^ {-1}) |G(x)|^2 \le  \bar \a_2 + \bar \b_2 |x|^2
\end{align*}
for all $a \in \R^d$.
\end{assumption}
When there is no super-linear term $G(x)$, we set $\theta = 0$ and $ \theta ^{-1} |G(x)|^2 =0$. Similarly, if linear term $G_1(x)$ is absent, we set $\theta  =  \infty $ and $ \theta  |G_1(x)|^2 =0$. Moreover, \eqref{linear_eq} implies
\begin{align*}
\la (2 x^T h(x) + |h(x)|^2) \le \la (|x|^2 + 2|h(x)|^2) \le 4\la K_1^2 + \la (4K_1 + 1)|x|^2, \qu  \forall x \in \R^d.
\end{align*}
Hence, by \As \ref{boundedness}, we have
\begin{align}\label{eq61}
2x^T f(x) + |g(x)|^2 + \la (2 x^T h(x) + |h(x)|^2) \le  \hat \a - \hat \b   |x|^2, \qu \forall x \in \R^d,
\end{align}
where $\hat \a  = \bar \a_1 + \bar \a_2 + 4 \la K_1^2  $ and $\hat \b  = \bar \b_1 - \bar \b_2 - \la (4 K_1^2 + 1)$.

\begin{theorem}\label{bounded_theorem}
Let \As \ref{Local_condition}, \ref{onesided_conditions}, \ref{monotone_condition} and \ref{boundedness} hold. Then for any initial value $x_0 \in \R^d$, the solution of the SDEs \eqref{eq0} satisfy
 \begin{align}\label{eq62}
 \mathop {\lim \sup }\limits_{t \to \infty } \E |x(t)|^2 \le \frac{\bar \a_1 + \bar \a_2 + 4 \la K_1^2 }{ \bar \b_1 - \bar \b_2 - \la (4 K_1^2 + 1)} .
 \end{align}
\end{theorem}
\pr Let $\tau_n$ , $\hat \a$, and $\hat \b $ be the same as before, $\sigma _n = t \we \tau_n $. For any $t \ge 0$, the It\^{o} formula gives that
\begin{align*}
\E \Big[ e^{\hat \b \sigma _n} |x( \sigma _n) |^2\Big ] &  = |x_0|^2 + \E \int_0^{\sigma _n }
 e^{\hat \b \sigma _n }  \Big (  2 x^T(s) f(x(s)) + |g(x(s))|^2   \no
& +   2 x^T(s) h(x(s)) + |h(x(s))|^2    + \hat \b | x(s)|^2\Big) ds.
\end{align*}
By \eqref{eq61}, we have
\begin{align*}\label{}
\E \Big[ e^{\hat \b\sigma _n } |x( \sigma _n) |^2\Big ] & \le  |x_0|^2 + \hat \a \int_0^{t} e^{ \hat \b s} ds  = |x_0|^2 + \frac{\hat \a}{\hat \b} (e^{ \hat \b t^-}  -1).
\end{align*}
Letting $n \to \infty$, we have
\begin{align*}\label{}
\E \Big[ e^{\hat \a t} |x( t) |^2\Big ] & \le  |x_0|^2 + \frac{\hat \a}{\hat \b} (e^{ \hat \b t^-}  -1)
\end{align*}
which implies
\begin{align*}\label{}
\E  |x( t) |^2\ \le \frac{|x_0|^2 }{e^{\hat \b t^-}} + \frac{\hat \a}{\hat \b} .
\end{align*}
Thus, the proof is complete. $\Box$

\begin{lemma}\label{lem_aa}
Let $0 < A <1$ and $B \ge0$. If
\begin{align}\label{aa1}
 D_k \le A D_{k-1} +B,\qu  \textrm{for} \qu   k=0,1,2, \cdots  .
\end{align}
 Then
\begin{align}\label{aa2}
 \mathop {\lim \sup }\limits_{k \to \infty } D_k \le \frac{B}{1-A}.
\end{align}
\end{lemma}
\pr The proof is given in the Appendix. $\Box$

\begin{theorem}\label{bounded_theorem2}
Let \As \ref{Local_condition}, \ref{onesided_conditions}, \ref{monotone_condition} and \ref{boundedness} hold. Then for any $\epsilon \in \big (0, \bar \beta_1 - \bar \b_2 - \max(\la (4K_1^2 +1), 2 \la K_1 (2 + K_1)) \big )$, there is a $\hat \DD \in (0, \Dst) $ and any initial value $x_0 \in \R^d$, the solution of the partially truncated EM method satisfy
 \begin{align}\label{eeee}
 \mathop {\lim \sup }\limits_{k \to \infty } \E |\XD(\tk)|^2 \le \frac{\bar \a_1 + \bar \a_2 + 2 \la K_1(2+K_1)+ \epsilon}{ \bar \b_1 - \bar \b _2 - 2 \la K_1(2+K_1) -\epsilon} .
 \end{align}
\end{theorem}
\pr Fix $\veps \in (0, \bar \b_1 - \bar \b_2)$. In the same way as Theorem 5.3 in \cite{Guo2017partial} was proved,  we have
\begin{align}\label{eq4022}
2x^T f_{\DD} (x) + |g _ {\DD}(x)|^2 \le \bar \a_1 + \bar \a_2  -( \bar \b _1 - \bar \b_2 - 0.5 \epsilon) |x|^2, \qu \forall x \in \R^d.
\end{align}
as long as $ \DD \in (0, \hat \DD _1)$, where $ \hat \DD _1 \in (0, \Dst)$ is sufficiently small for which
\begin{align}\label{eq51}
\frac{\bar \a_2}{ (\mu^{-1}(\ph(\hat \DD_1 )))^2} \le 0.5 \epsilon.
\end{align}
Using the independence of the Poisson increments and \eqref{linear_eq} as well as Lemma \ref{moment of xdt}, we obtain that
\begin{align}\label{eq503}
\E |h (\xdtk) \DD N_k|^2 & \le 2  K_1^2 \E(1+ | \xdtk|^2) \E | \DD N_k|^2 \no
& \le   2  K_1^2 \la \DD (1+ \la \DD ) \E (1+| \xdtk|^2 )\no
& \le  2  K_1^2 \la \DD \E | \xdtk|^2 + 2K_1^2 \la \DD + C \DD^2,
\end{align}
\begin{align}\label{eq504}
  2 \E | \xdtk  h (\xdtk)\DD N_k | & \le 2K_1 \E (|\xdtk|(1+ |\xdtk|)) \E | \DD N_k| \no
  & \le 4 K_1  \la \DD  \E (1+|\xdtk|^2) \no
  & \le  4 K_1  \la \DD \E |\xdtk|^2 +  4 K_1  \la \DD
\end{align}
and
\begin{align}\label{eq505}
2 \E |\DD \fD  (\xdtk) h (\xdtk)\DD N_k |  & \le 2 K_1 \DD \E ( (1+ |\xdtk |)|\fD (\xdtk) | ) \E |\DD N_k| \no
 & \le K_1 \la \DD ^2 (  \E (1+| \xdtk|^2) + \E | \fD (\xdtk)| ^2 )  \no
 & \le   K_1 \la \DD ^2 \E | \fD (\xdtk)| ^2 +C \DD ^2.
\end{align}
Fix $x_0 \in \R^d$ arbitrarily. For any $ \DD \in (0, \hat \DD_1)$,
substituting \eqref{eq503}-\eqref{eq505} into \eqref{eq406} gives
\begin{align}\label{eq508}
\E | \XD ( \tkk) |^2 &  \le  \E \Big ( |\xdtk|^2 +  2 \xdtk ^T  \fD (\xdtk)  \DD  + |\gD (\xdtk) |^2   \DD  \Big )  \no
& \qu + 2 \la K_1 (2 + K_1) \DD  \E | \xdtk |^2  +  (1+ K_1 \la )  \DD ^2  \E |\fD (\xdtk)| ^2 \no
& \qu + 2 \la K_1 (2 + K_1) \DD  + C \DD ^2  \no
&  \le     (1- (\bar \b_1  - \bar \b_2 -  2 \la K_1 (2 + K_1) - 0.5 \veps) \DD )   \E |\xdtk|^2  \no
                    & \qu + ( \bar \a_1 + \bar \a _2 + 2 \la K_1 (2 + K_1)  ) \DD  + C \DD^2 \no
                    & \qu + (1+ K_1 \la )  \DD ^2  \E |\fD (\xdtk)| ^2,
\end{align}
where \eqref{eq4022} have been used.
By \eqref{linear_eq} and \eqref{HD}, we have
\begin{align*}
|f_{\DD}(x)|^2 \le 2 |F_1(x)|^2 + 2 |F(x)|^2 \le 4K_1^2 (1+ |x|^2) + 2 (\ph(\DD))^2, \qu \forall x \in \R^d.
\end{align*}
Hence, by \eqref{delta_h_relations}, we get
\begin{align*}
\DD|f_{\DD}(x)|^2  \le 4K_1^2 \DD(1+ |x|^2) + 2\DD ^ { 1/2 \we (\bar p -2) / \bar p
 }, \qu \forall x \in \R^d.
\end{align*}
Consequently, there is a $\hat \DD \in (0, \hat \DD_1)$ sufficiently small such that for any $\DD \in (0, \hat \DD)$, $\DD (\hat \b_1 - \hat\b_2 - \epsilon) <1$ and
\begin{align} \label{eq510}
C\DD + (1+K_1 \la)\DD|f_{\DD}(\XD(\tk))|^2 \le \veps + 0.5 \epsilon |\XD(\tk)|^2.
\end{align}
Thus, fix any $\DD\in (0, \hat \DD)$. Inserting \eqref{eq510} into \eqref{eq508} yields
\begin{align}\label{eq511}
\E | \XD ( \tkk) |^2 &  \le     (1- (\bar \b_1  - \bar \b_2 -  2 \la K_1 (2 + K_1) - \epsilon) \DD )   \E |\xdtk|^2  \no
                    & \qu + ( \bar \a_1 + \bar \a _2 + 2 \la K_1 (2 + K_1) + \epsilon ) \DD.
\end{align}
Applying Lemma \ref{lem_aa} to \eqref{eq511} gives the required assertion \eqref{eeee}.  $\Box$

\section{Examples}

\begin{example}\label{example1}
Consider the following scalar SDEs with jumps
 \begin{align}\label{ex2}
dx(t^-) =  - x^5(t^-)dt + x^2(t^-) dB(t) + x^2(t^-) d N(t),
 \end{align}
with the initial value $x(0) = 1 $,
where $B(t)$ is a scalar Brownian motion and $N(t)$ is a scalar Poisson process with jump intensity $\la = 0.5$. Obviously, we have
\begin{align*}
 f(x) = - x^5 , \qu g(x) =x^2, \qu h(x) = x^2.
\end{align*}
It is easy to check that \As \ref{2Polynomial_condition} is satisfied with
\begin{align} \label{exeq5}
|f(x) - f(y)| \ve  |h(x) - h(y)| \le \bar L_2 (1 + x ^{4 } + y^{4 } )|x -y|
\end{align}
for a constant $ \bar L_2  >0 $.
The elementary inequality gives
\begin{align}\label{exeq1}
& 2 \la (x -y ) (h(x) - h(y)) + \la |h(x)-h(y)|^2   \no
& = (x-y) (x^2 -y^2)  +  0.5|x^2 - y^2 |^2  \no
& \le  |x +y| |x -y|^2  + 0.5 |x +y|^2 |x-y|^2  \no
& \le   ( x^2 + y^2 +1 ) |x -y|^2 +   (x ^2 + y^2) |x-y|^2  \no
& = \B ( 2 ( x^2 + y^2)  +1 \B ) |x -y|^2.
\end{align}
Also, we have
\begin{align}\label{exeq2}
- (x^3y + x y^3 + x^2 y^2 ) = -xy (x^2 +y^2) - x^2 y^2  \le 0.5 (x ^4 + y^4).
\end{align}
By \eqref{exeq1} and \eqref{exeq2}, we have

\begin{align}\label{exeq3}
& 2(x-y ) (f(x) - f(y)) + |g(x) - g(y)|^2 +2 \la (x-y) |h(x) - h(y)| + \la |h(x) - h(y)|^2 \no
& \le (x-y ) \B ( - (x -y)(x^4 + x^3y + x^2 y^2 + x y ^3 + y ^4 )  \B)  + | x^2 - y^2 |^2 \no
& \qu  + \B ( 2 ( x^2 + y^2)  +1 \B ) |x -y|^2 \no
& \le - 0.5 (x^4 + y^4) |x -y|^2 + 2(x ^2 +y^2)|x -y|^2 + \B ( 2 ( x^2 + y^2)  +1 \B ) |x -y|^2\no
& =  \B ( - 0.5 (x^4 + y^4) +  4 ( x^2 + y^2)  +1 \B ) |x -y|^2\no
& \le 17 |x -y|^2.
\end{align}
 Note that in the last inequality the elementary inequality $$k a b \le \frac{1}{2} a^2 + \frac{1}{2} k^2 b ^2,  \qu \forall  a,b , k \in \R $$ has been used. Hence, \As \ref{2onesided_conditions} is satisfied.
For \As \ref{2monotone_condition}, we have
\begin{align}\label{exeq4}
 & 2 x^ T f(x)   +  | g(x)|  ^2  + \la ( 2 x^ T h(x) + |h(x)|^2)  \no
 & = -2x^6 + x^4 + 0.5( 2x^3 + x^4) \no
 & = (-x^6 + \frac{3}{2} x^4 )+  ( x^3   - x^6  )  \no
 & = x^2 \B ( -(x^2 - \frac{3}{4})^2  + \frac{9}{16} \B ) +  ( x^3   - x^6  )    \no
 & \le \frac{9}{16}x^2 + \frac{1}{4},
\end{align}
where we use the elementary inequality $ ab \le a^2 + b^2 /4 $ in the last inequality.
By \eqref{exeq5}, we can choose $\mu(n) = n^5$ such that
$$\sup _{ |x| \le n} (|f(x)|\ve |g(x)| \ve |h(x)|) = \sup _{ |x| \le n} (|x^5| \ve |x^2| \ve |x^2|) \le n^5 ,\qu \forall n \ge 1.  $$
Letting $\bar \g =4$, $r =2/(2 + \bar \g) = 1/3$ and $ \veps = 1/4$.  If we choose $\ph ( \DD) =  \DD ^ {-\veps} =  \DD ^ { - 1/4}$, then all the conditions in \eqref{2delta_h_relations} hold for all $\Dst \in (0, 1].$
Hence, the truncating factor is defined as $ \mu ^ {-1} ( \ph (\DD)) = \DD ^ { - 1/20}$ and the truncated functions are defined as
\begin{align*}
f_{\DD} (x) = f( ( |x | \we  \DD ^ { - 0.05}) \frac{x}{|x|}) \no
 g_{\DD} (x) = g( ( |x | \we  \DD ^ { - 0.05}) \frac{x}{|x|}) \no
 h_{\DD} (x) = h( ( |x | \we  \DD ^ { - 0.05}) \frac{x}{|x|}).
\end{align*}
For the given step size $\DD$ and the time $T $, the $ X _ {k+1}$ is calculated by
\begin{align}\label{Tem}
X _ {k+1}  = X _ {k} + f_{\DD} (X _ {k} )  \DD  +g_{\DD} (X _ {k} ) \DD B_k + h_{\DD} (X _ {k} ) \DD N_k, \quad 0 \le k \le T/\DD -1,
\end{align}
with $X _ 0 =1$.
Then, by Corollary \ref{cor2}, the truncated EM scheme converges strongly with rate $  r/2 - r \veps = 1/(4 + 2 \bar \g) = 1/12 $ to the true solution of the SDE \ref{ex2}.

As SDEs \ref{ex2} does not have any explicit solutions, the scheme \eqref{Tem} with step size $2^{-16}$ is treated as the solution of the SDEs \ref{ex2} in the numerical experiment. The number of simulation is $500$.  The $ \mathcal L ^ {1/3}$ errors at time $T =4 $ ,that is $$ \B( \E \big |X(T) - x_{\DD}(T) \big |^ {1/3}  \B)^3   \approx \B ( \frac{1}{500} \sum _ {i=1}^ {500} \big | [X(T)]^i - [x_{\DD}(T)]^i \big |^{1/3} \B )^3,$$
with step sizes $2^{-15}$, $2^{-14}$, $2^{-13}$, $2^{-12}$ and $2^{-11}$ are displayed in Fig \ref{ex1_fig} .
The numerical result shows that our numerical findings are consistent with the theoretical results obtained in this paper.
\begin{figure}[!t]
\centering
  \includegraphics[width=8cm,height=6cm]{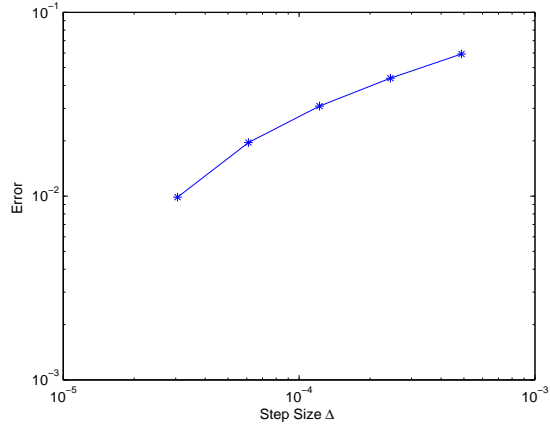}
\caption{$\mathcal L ^ {1/3} $-convergence of truncated EM scheme \ref{Tem} of SDEs \ref{ex2} }
\label{ex1_fig}
\end{figure}
\end{example}

\begin{figure}[!t]
\centering
\begin{minipage}[h]{0.45\textwidth}
   \centering
   \includegraphics[width=6cm,height=5cm]{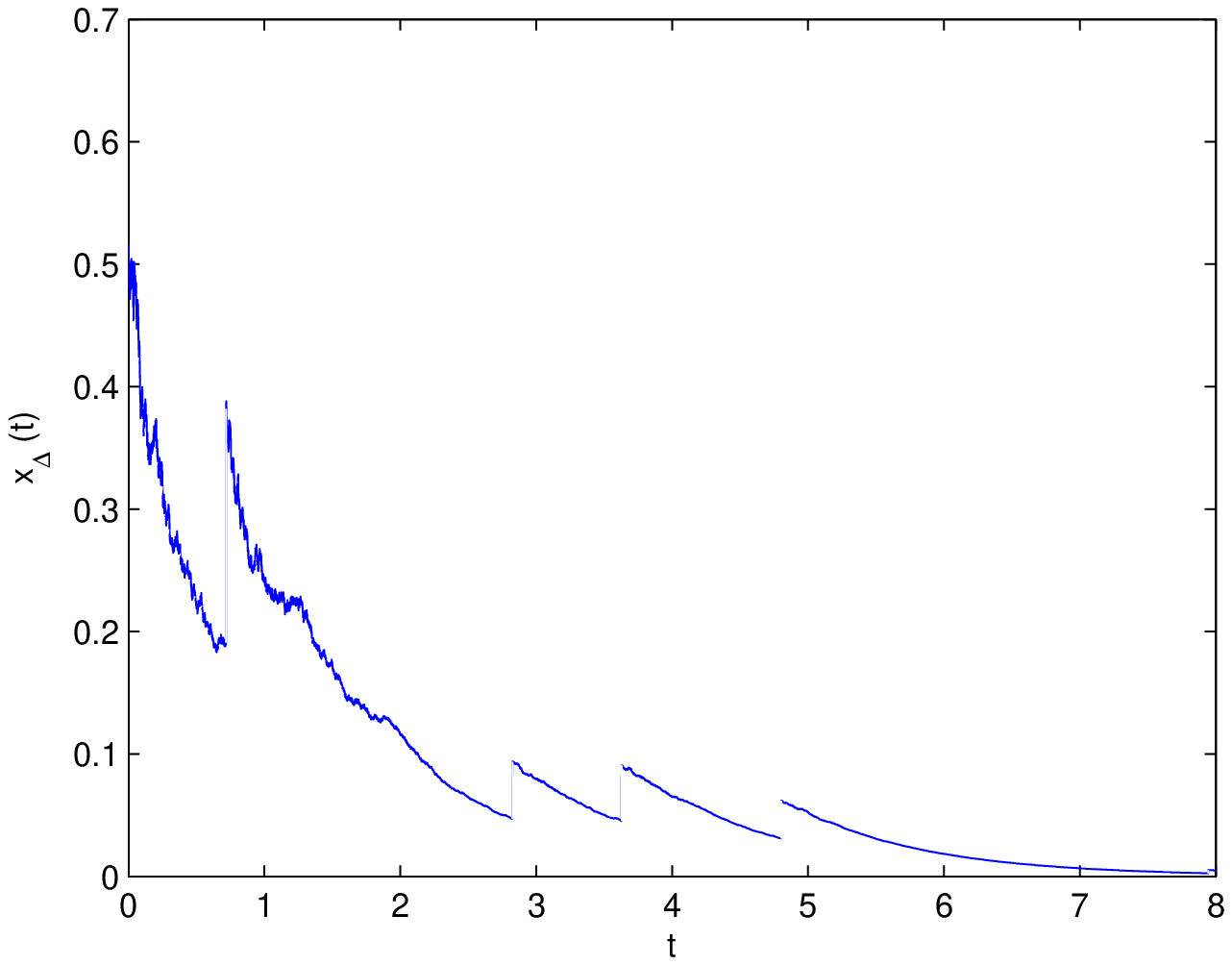}
\caption{Simulation of one path in Example 5.2}
\label{fig3}
\end{minipage}
\hfil
\begin{minipage}[h]{0.45\textwidth}
\centering
  \includegraphics[width=6cm,height=5cm]{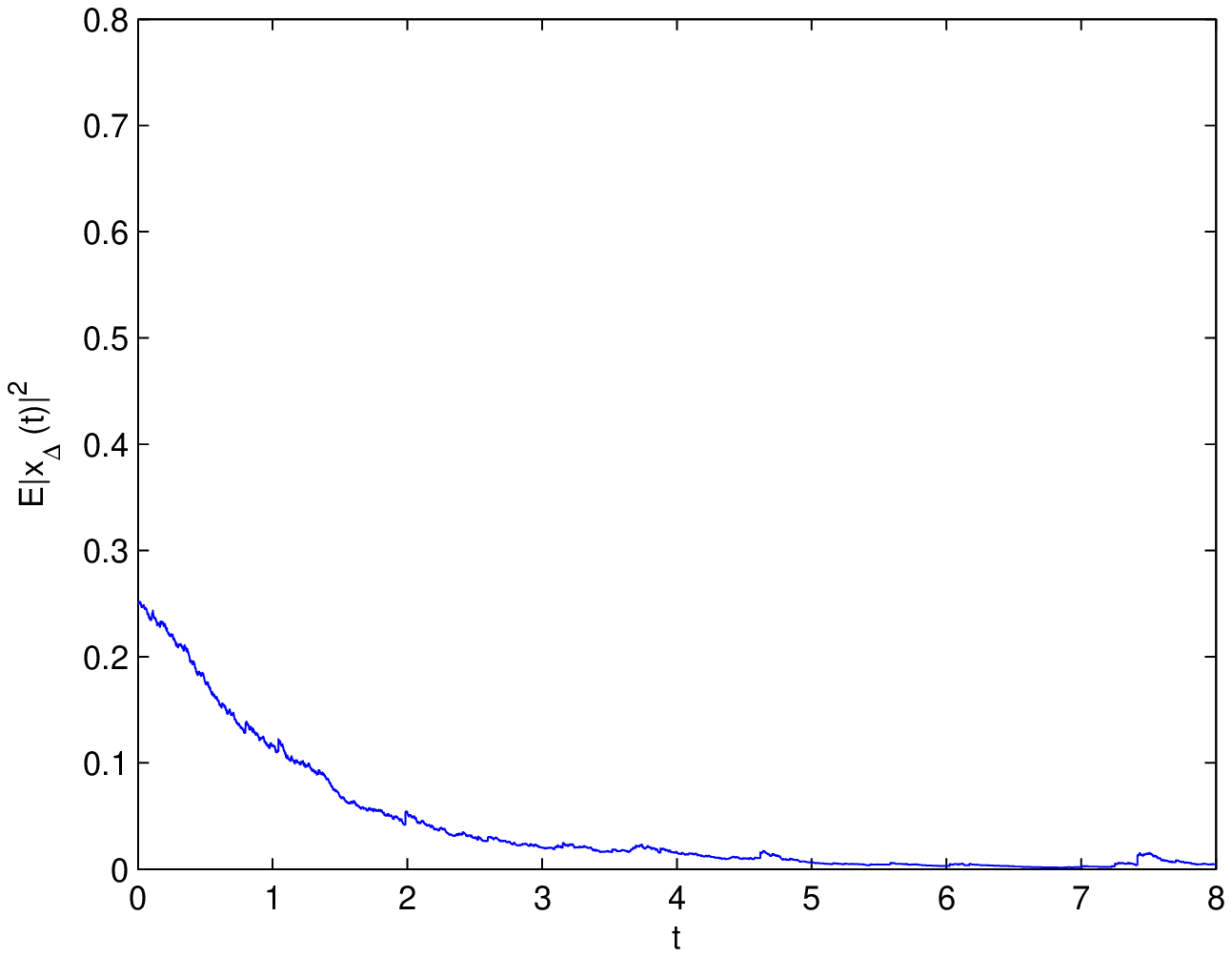}
\caption{Mean square of 1000 paths in Example 5.2  }
\label{fig4}
\end{minipage}
\end{figure}

\begin{example}\label{example2}
Consider the following scalar SDEs with jumps
 \begin{align}\label{ex3}
dx(t^-) =  - (x(t^-) +x^5(t^-))dt + x^2(t^-) dB(t) + x(t^-) d N(t),
 \end{align}
with the initial value $x(0) = 0.5 $,
where $B(t)$ is a scalar Brownian motion and $N(t)$ is a scalar Poisson process with jump intensity $\la = 0.5$. Obviously, we have
\begin{align*}
F_1(x) = -x, \qu F(x) = - x^5 , \qu G_1(x) =0, \qu G(x) = x^2, \qu h(x) = x,
\end{align*}
and
\begin{align*}
|F_1(x)| \ve |G_1(x)| \ve |h(x)|  =|x|, \qu \textrm{with} \qu K_1=1.
\end{align*}
 Setting $ \theta = \infty$ gives
$$ 2 x F_1(x) + ( 1+ \theta ) |G_1(x)|^2  = -2 x^2,$$
and
$$ 2 x F(x) + (1 + \theta ^ {-1} )|G(x)|^2 = -2 x^6 + x^4 \le  -2 x^2 \B ( x^2 - \frac{1}{4} \B )^2 + \frac{1}{8} x^2 \le \frac{1}{8} x^2. $$
Hence, \As \ref{stability_condition} is satisfied with $ \a _1 =2$ and $ \a _2 = 1/8$.
Moreover, we have
\begin{align*}
(x-y) (F(x) - F(y)) + \frac{ \bar r -1}{2} |G(x) - G(y)|  \le \B ( 1 +  \frac{(\bar r -1)}{4} \B) |x-y|^2, \
\qu \forall x  \in \R
\end{align*}
which means that  \As \ref{onesided_conditions} is satisfied for any $\bar r$.
Also, we can check that \As \ref{monotone_condition} is fulfilled for any $\bar p$ (see \cite{Guo2017partial}). By Theorem \ref{stability_theorem}, the SDE \ref{ex3} is stable exponentially in the mean square sense for any initial value $x_0 \in \R$, and the solution of SDE \ref{ex3} satisfies
\begin{align*}
\E |x(t^-)|^2 \le |x_0|^2 e^{-(\a_1 - \a_2 - \la K_1 (2 + K_1))t^-}   =  |x_0|^2 e ^{-0.375 t^- }, \qu \forall t \ge 0.
\end{align*}
Letting  $r =2$, $\bar r =3$. If we choose $ \mu (n ) = n^5$,$\g =4,  \ph (\DD) = \DD ^ {- 1/ 40}$, and $ p = 40$, then by Corollary \ref{cor1}, the numerical solution will converge strongly to the true solution in $\mathcal L^2$ with convergence rate $ 1/( 1 + \g) - r/p  = 0.15 $.
Finally, by Theorem \ref{stability_theorem2}, for any $ \epsilon \in (0, 0.375)$, there exists a $\hat \DD
\in (0, \Dst) $ such that for all $ \DD \in (0, \hat \DD)$ and any initial value $x_0 \in \R ^d$, the solution of the truncated EM method \eqref{EM scheme} satisfies
\begin{align*}
\E |\XD (\tk)|^2 \le |x_0|^2 e ^{-(\a_1 - \a_2 - \la K_1 (2 + K_1) - \epsilon) \tk  } = |x_0|^2 e ^{-(0.375 - \epsilon) \tk  }, \qu \forall  k \ge 0.
\end{align*}
Fig \ref{fig3} and Fig \ref{fig4} demonstrate the mean square exponential stability of the truncated EM method.
\end{example}

\begin{figure}[!t]
\centering
\begin{minipage}[h]{0.45\textwidth}
   \centering
   \includegraphics[width=6cm,height=5cm]{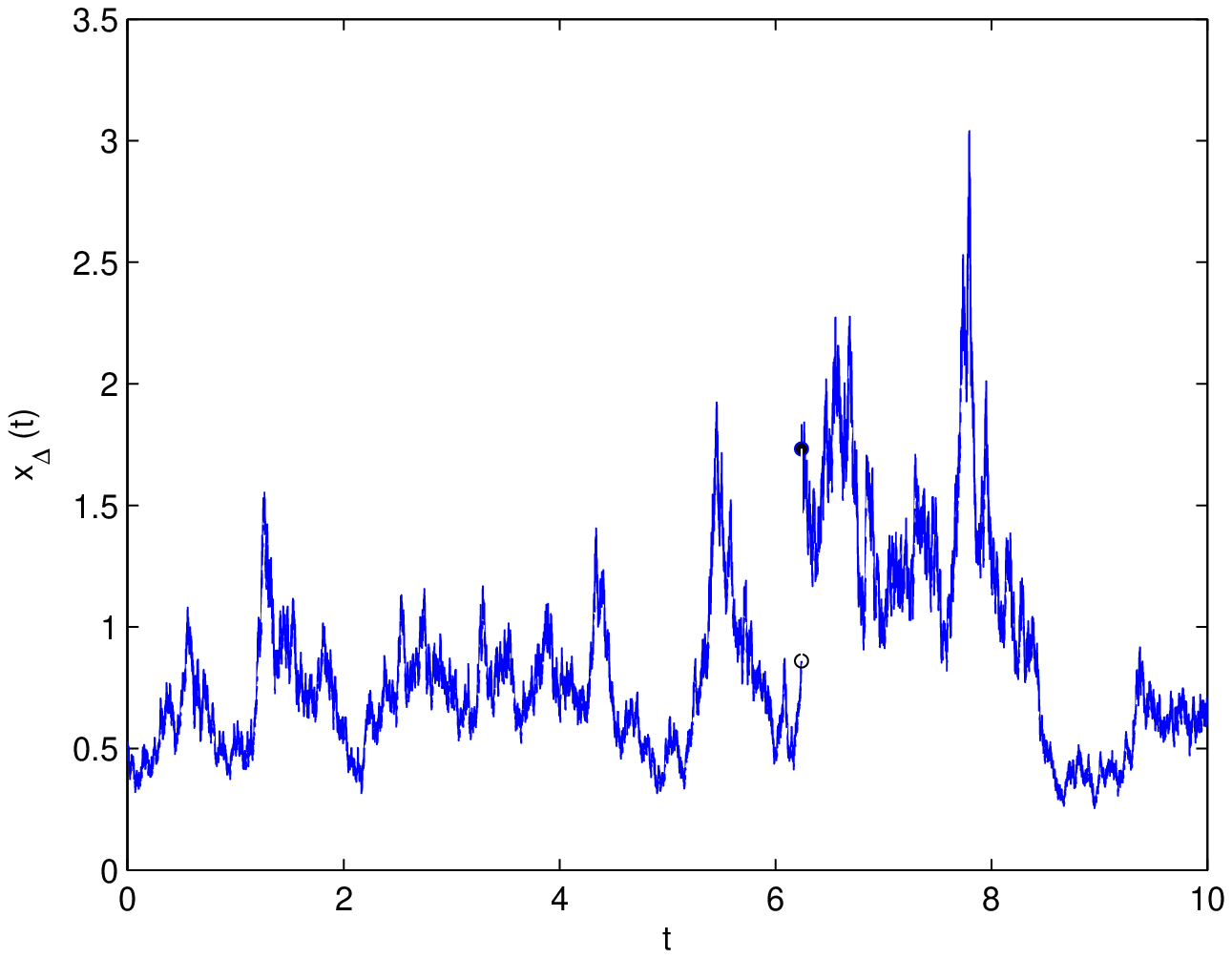}
\caption{
Simulation of one path in Example 5.3}
\label{fig5}
\end{minipage}
\hfil
\begin{minipage}[h]{0.45\textwidth}
\centering
  \includegraphics[width=6cm,height=5cm]{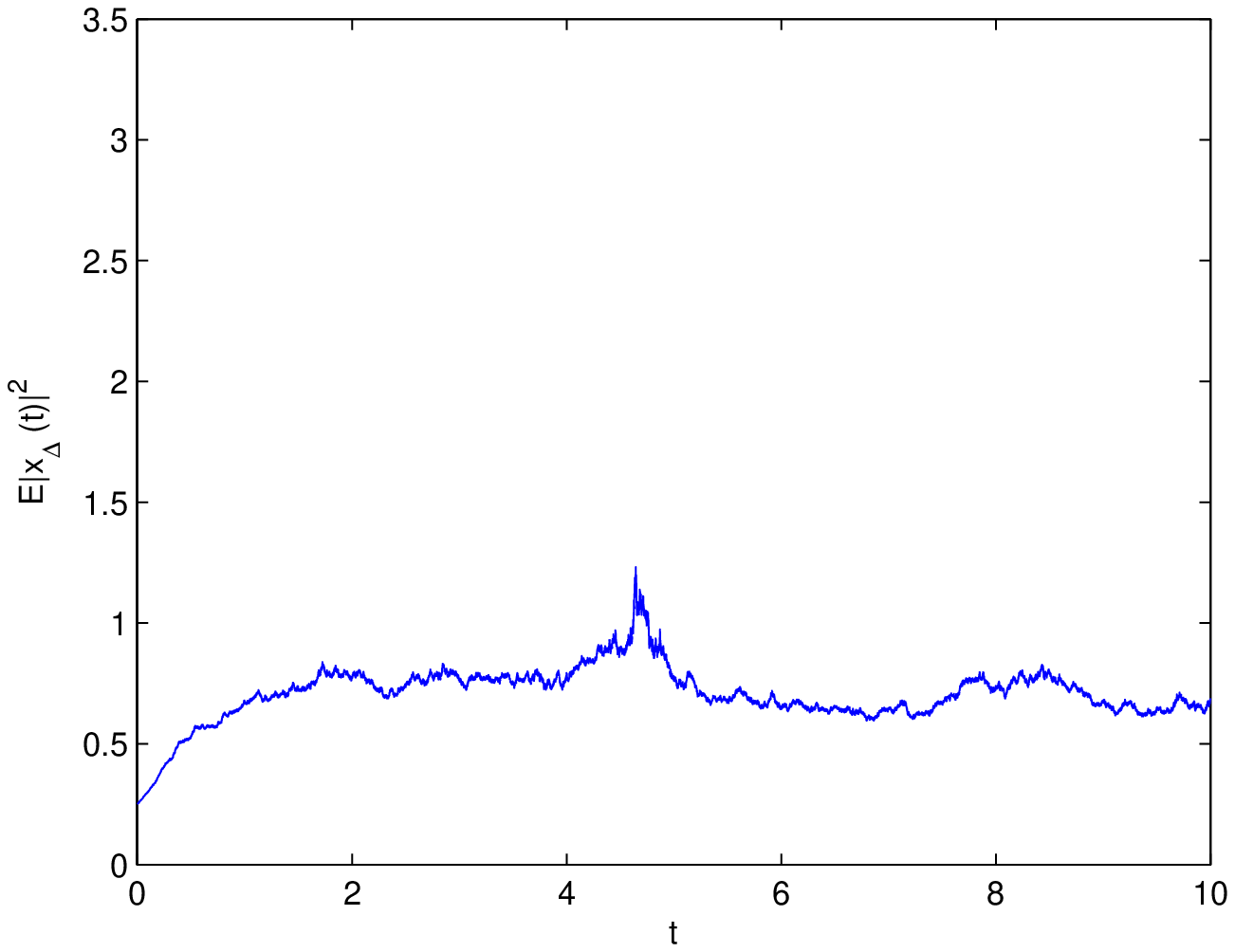}
\caption{Mean square of 1000 paths in Example 5.3 }
\label{fig6}
\end{minipage}
\end{figure}

\begin{example}
Consider the following scalar SDEs with jumps
 \begin{align}\label{ex4}
dx(t^-) =   (x(t^-) - x^3(t^-))dt + x(t^-) dB(t) + x(t^-) d N(t),
 \end{align}
with the initial value $x(0) = x_0 \in \R$,
where $B(t)$ is a scalar Brownian motion and $N(t)$ is a scalar Poisson process with jump intensity $\la = 0.1$.
We decompose the drift and diffusion coefficient in the form with
\begin{align} \label{3eq1}
F_1(x) = -2x, \qu F(x) = 3x- x^3 , \qu G_1(x) =x, \qu G(x) = 0, \qu h(x) = x,
\end{align}
which means
\begin{align*}
|F_1(x)| \ve |G_1(x)| \ve |h(x)|  = 2|x|, \qu \textrm{with} \qu K_1=2.
\end{align*}
Setting $ \theta = 0 $ gives
$$ 2 x F_1(x) + ( 1+ \theta ) |G_1(x)|^2  = -3 x^2,$$
and
$$ 2 x F(x) + (1 + \theta ^ {-1} )|G(x)|^2 = 2x(3x - x^3) = -2 ( x^2 - 1.5)^2 + 4.5 \le 4.5 . $$
Hence, \As \ref{boundedness} is satisfied with
\begin{align}\label{cond111}
\bar \a _1 = 0 , \qu  \bar \b _1 = 3, \qu   \bar \a_2 = 4.5 ,\qu \textrm{and} \qu  \bar \b _2 = 0.
\end{align}
It is easy to check that coefficients of the SDE \ref{ex4} with their decompositions in \eqref{3eq1} satisfy \As \ref{Local_condition}, \ref{onesided_conditions} and \ref{monotone_condition} for any $ \bar p >2$.
Applying Theorem \ref{bounded_theorem} gives that for any initial value $x_0 \in \R^d$, the solution of SDE \ref{ex4} satisfies
 \begin{align}\label{3eq222}
 \mathop {\lim \sup }\limits_{t \to \infty } \E |x(t)|^2 \le \frac{\bar \a_1 + \bar \a_2 + 4 \la K_1^2 }{ \bar \b_1 - \bar \b_2 - \la (4 K_1^2 + 1)}  \approx 4.69.
 \end{align}
Moreover, we can choose $ \mu (n) = 4 n^3$ and $ \ph(\DD) = \DD ^ { - \veps}$ and let $ r =2$, $ \bar \g =2$, $ p = 50$ as well as $ \veps = 1/50$ to define the
numerical solution $ X_ {\DD} ( t_k )$ by the partially truncated EM method. By Theorem \ref{lemma3.2}, this solution of truncated EM will converge to the true solution in $\mathcal L ^ 2$ with convergence rate $ 1/(1 + \bar \g) - r/p \approx 0.2933 $. Finally, by Theorem \ref{bounded_theorem2},
%
for any $ \epsilon \in (0, 1.3)$, there exists a $\hat \DD
\in (0, \Dst) $ such that for all $ \DD \in (0, \hat \DD)$ and any initial value $x_0 \in \R ^d$, the numerical solution satisfies
 \begin{align*}
 \mathop {\lim \sup }\limits_{k \to \infty } \E |\XD(\tk)|^2 \le \frac{\bar \a_1 + \bar \a_2 + 2 \la K_1(2+K_1)+ \epsilon}{ \bar \b_1 - \bar \b _2 - 2 \la K_1(2+K_1) -\epsilon}  = \frac{6.1 + \epsilon}{ 1.4 - \epsilon}.
 \end{align*}
The asymptotic boundedness of the numerical method is shown in Fig \ref{fig5} and Fig \ref{fig6}.
\end{example}

\section{Conclusions and future research} \label{secon}

In this paper, the truncated EM method is investigated for SDEs driven by both Brownian motions and Possion jumps. Both the finite time convergence and asymptotic behaviours of the method are studied. The $\mathcal L^r(r \geq 2)$ strong convergence is proved when the drift and diffusion coefficients satisfy super-linear growth condition and the coefficient for Possion jumps satisfies linear growth condition. When $0 < r <2 $, we are able to prove the $ \mathcal L^r$ convergence of the methods to SDEs with all the three coefficients allowing to grow super-linearly.
\par
In the future works, we will report on the SDEs driven by L\'evy process and the $\mathcal L^r$ convergence for SDEs whose all the three coefficients can grow super-linearly.

\begin{appendix}
\section{Proof of Lemma \ref{lemma303}}
\pr By the \Ito and \eqref{eq28}, we have
\begin{align*}\label{}
\E | x( t \we \tr)|^2 &  \le \xx + \E  \int _0 ^ { t \we \tr}  K_3 (1+ |x(s)|^2) ds  \no
   & \qu + \lambda \E \int _0 ^ { t \we \tr} ( 2 x(s)^T h(x(s)) + |h(x(s))|^2  ) ds \no
   & \le \xx +( K_3 + 2 \la (2K_1 + K_1^2)) \intt  \E (1+| x(s \we \tr )|^2 )  ds
\end{align*}
for any $ 0 < t < T$. The Gronwall inequality shows \begin{align*}
         \E |x(T \we \tr ) |^2  \le C.
                                                    \end{align*}
This implies
$$ \PP (\tr \le T) \le \frac{C}{n^2} .$$  Thus, the proof is complete. $\Box$
\section{Proof of Lemma \ref{lemma34}}
\pr We write $ \rho _ {\DD , n} = \rho $ for simplicity. For $0 \le t \le T$, the It\^{o} formula gives
\begin{align}\label{eq123}
\E | x _ {\DD} (t \we \rho )|^2   & = \xx + \E  \int _0 ^ {t \we \rho} \Big (  2 x _{\DD} ^ T (s) \fdxs + |\gdxs|^2  \Big )ds \no
& \qu +  \lambda \E  \int _0 ^ {t \we \rho} \Big (  2 x _{\DD} ^ T (s)  h(\bxds) + | h(\bxds) |^2  \Big ) ds \no
& = \xx + \E  \int _0 ^ {t \we \rho} \Big (  2\bar x _{\DD} ^ T (s) \fdxs + |\gdxs|^2  \Big )ds \no
& \qu + \E  \int _0 ^ {t \we \rho}  2 ( \xds - \bxds)^ T \fdxs   ds \no
& \qu +  \lambda \E  \int _0 ^ {t \we \rho} \Big (  2 x _{\DD} ^ T (s)   h(\bxds) + | h(\bxds) |^2  \Big ) ds.
\end{align}
By \eqref{linear_eq}, we obtain
\begin{align}\label{aa}
& \E  \int _0 ^ {t \we \rho} \Big (  2 x _{\DD} ^ T (s)  h(\bxds) + |h(\bxds) |^2  \Big ) ds \no
  & \le \E  \int _0 ^ {t \we \rho} \Big ( |\xds |^2+ 2 |h(\bxds) |^2  \Big ) ds  \no
  & \le \E  \int _0 ^ {t \we \rho} \Big ( |\xds |^2+ 4K_1^2 (1+ | \bxds |^2) \Big ) ds  \no
    & \le 4K_1^2T + (8K_1^2T+1) \E  \int _0 ^ {t \we \rho} |\xds |^2 ds +8K_1^2T \E  \int _0 ^ {t \we \rho} |\xds - \bxds |^2 ds .
\end{align}
Substituting this into \eqref{eq123} and applying \ref{eq209}, we have
\begin{align}\label{bb}
\E | x _ {\DD} (t \we \rho )|^2  & \le \xx + \int _0 ^ {t \we \rho} 2K_4 (1 + |\bxds|^2)ds  + 4 \la K_1^2 T \no
 & \qu + \E  \int _0 ^ {t \we \rho}  2 ( \xds - \bxds)^ T \fdxs   ds  \no
 & \qu +  \la (8K_1^2T+1) \E  \int _0 ^ {t \we \rho} |\xds |^2 ds + \la 8K_1^2T \E  \int _0 ^ {t \we \rho} |\xds - \bxds |^2 ds \no
 & \le ( \xx + 2K_4T + 4 \la L_1^2T )+ (4K_4 + \la (8K_1^2T+1)) \intt \E | x _ {\DD}(s \we \rho )|^2 ds \no
 & \qu + (4K_4 + 8\la K_1^2T) \int_0^T \E | \xds - \bxds|^2 ds \no
 & \qu + 2 \E \int_0 ^{t \we \rho} |\xds - \bxds| |\fdxs| ds.
\end{align}
By Lemma \ref{lem17}, we have
\begin{align*}
 \int_0^T \E |\xds - \bxds|^2 ds \le C.
\end{align*}
By \eqref{linear_eq}, we have
\begin{align}\label{}
& \E \int_0 ^{t \we \rho} |\xds - \bxds| |\fdxs| ds \no
& \le K_1 \E \int_0 ^{t \we \rho} |\xds - \bxds| (1+|\bxds|) ds +I_5 \no
& \le  C \Big ( \E  \int _0 ^ {t \we \rho} |\xds - \bxds |^2 ds + \intt \E | x _ {\DD}(s \we \rho )|^2 ds +1  \Big ) + I_5
\end{align}
where \begin{align*}
I_5 = \E \int_0 ^T |\xds - \bxds| |F_{\DD} (\bxds)| ds .
      \end{align*}
Using Lemma \eqref{lemma33}, condition \ref{delta_h_relations} and \eqref{HD} gives
\begin{align*}
I_5 &  \le  \ph(\DD)  \int _0 ^ T \Big ( \E |\xds- \bxds|^2 \Big )^ {1/2} ds \no
    & \le  C (\ph(\DD))^2 \DD ^ {1/2} =   C (\ph(\DD) \DD ^ {1/4} )^2 \le C.
\end{align*}
Hence, we have
\begin{align*}
 \E  |x _ {\DD} (t \we \rho )|^2  \le C \Big (1 +  \intt \E  |x _ {\DD} (s \we \rho )|^2  ds \Big  ).
\end{align*}
The Gronwall inequality gives
\begin{align*}
\E |x _ {\DD} (T \we \rho )|^2 \le C,
\end{align*}
which implies \eqref{eq35}. Thus, the proof is complete. $\Box$

\section{Proof of Lemma \ref{2lem 24}}
\pr Fix any $\DD \in (0, \Dst ]$. Recall that $\mu ^{-1} (\ph(\Dst))\ge 1$, we get  $\mu ^{-1} (\ph(\DD))\ge 1$.
For $x \in \R^d$ with $|x| \le  \mu ^{-1} (\ph(\DD))$, by the definition of the truncated function, we obtain the required
assertion \eqref{ddd}.  For $x \in \R^d$ with $|x| >  \mu ^{-1} (\ph(\DD))$, by \eqref{2monotone_condition}, we have
\begin{align*}
& \qu  2  x^{T} f_ {\Delta} (x) +  |g_{\Delta}(x)|^2   +  \la ( 2 x^T h_{\DD}(x) + |h_{\DD} (x)|^2 )\no
& = 2 (x - \pi _{\DD}(x)) ^T  f_ {\Delta} (x) + 2 \la (x - \pi _{\DD}(x)) ^T  h_ {\Delta} (x) \no
& \qu + 2 \pi _{\DD}(x)^T f_ {\Delta} (x) + |g_{\Delta}(x)|^2+ 2 \la \pi _{\DD}(x)^T h_{\DD}(x) + \la |h_{\DD}(x)|^2 \no
& \le \B ( \frac{|x|}{\mu ^{-1} (\ph(\DD))} -1  \B ) \B (  2 \pi _{\DD}(x)^T f(\pi _{\DD}(x)) + 2 \la  \pi _{\DD}(x)^T h(\pi _{\DD}(x)) \B)  + \bar K (1+ |\pi _{\DD}(x)|^2) \no
& \le  \B ( \frac{|x|}{\mu ^{-1} (\ph(\DD))} -1  \B ) ( \bar K (1+ |\pi _{\DD}(x)|^2) ) + \bar K (1+ |\pi _{\DD}(x)|^2) \no
& = \frac{|x|}{\mu ^{-1} (\ph(\DD))} \bar K ( 1+ |\mu ^{-1} (\ph(\DD))|^2 )  \no
& =  \bar K |x|  \B ( \frac{1 }{\mu ^{-1}  (\ph(\DD))}+ |\mu ^{-1} (\ph(\DD))| \B )  \no
& \le \bar K |x| (1 + |x|)  \le  \bar K  (1 + |x|) ^ 2 \le 2 \bar K (1 + |x|).
\end{align*}
Thus, we complete the proof. $\Box$

\section{Proof of Lemma \ref{lem_aa}}
\pr \eqref{aa1} is equivalent to the following expression
\begin{align*}
 D_k + \frac{B}{A-1} \le A \left( D_{k-1} +\frac{B}{A-1} \right ),\qu  \textrm{for} \qu   k=0,1,2, \cdots  .
\end{align*}
Hence, we have
\begin{align*}
 D_k + \frac{B}{A-1} \le A^k \l ( D_{0} +\frac{B}{A-1} \r ).
\end{align*}
It follows
\begin{align*}
 D_k   \le A^k \l ( D_{0} +\frac{B}{A-1} \r ) + \frac{B}{1-A}.
\end{align*}
Recalling $0<A<1$ and taking $k \to \infty$, we obtain the required assertion \ref{aa2}.
Thus, the proof is complete.  $\Box$
\end{appendix}
\section*{Acknowledgment}

This work was supported in part by the Natural Science Foundation of China (No. 71571001, 61703003).

\section*{References}
\bibliographystyle{model1-num-names}
\bibliography{refs}

\begin{thebibliography}{25}
\expandafter\ifx\csname natexlab\endcsname\relax\def\natexlab#1{#1}\fi
\providecommand{\bibinfo}[2]{#2}
\ifx\xfnm\relax \def\xfnm[#1]{\unskip,\space#1}\fi
\bibitem[{Allen(2007)}]{ALL2007a}
\bibinfo{author}{E.~Allen}, \bibinfo{title}{Modeling with {I}t\^o Stochastic
  Differential Equations}, \bibinfo{publisher}{Springer},
  \bibinfo{address}{Dordrecht}, \bibinfo{year}{2007}.
\bibitem[{Mao(2008)}]{Mao2007book}
\bibinfo{author}{X.~Mao}, \bibinfo{title}{Stochastic Differential Equations and
  Applications}, \bibinfo{publisher}{Horwood}, \bibinfo{edition}{2nd} edition,
  \bibinfo{year}{2008}.
\bibitem[{Higham and Kloeden(2005)}]{HK2005a}
\bibinfo{author}{D.~J. Higham}, \bibinfo{author}{P.~E. Kloeden},
\newblock \bibinfo{title}{Numerical methods for nonlinear stochastic
  differential equations with jumps},
\newblock \bibinfo{journal}{Numer. Math.} \bibinfo{volume}{101}
  (\bibinfo{year}{2005}) \bibinfo{pages}{101--119}.
\bibitem[{Higham and Kloeden(2006)}]{HK2006a}
\bibinfo{author}{D.~J. Higham}, \bibinfo{author}{P.~E. Kloeden},
\newblock \bibinfo{title}{Convergence and stability of implicit methods for
  jump-diffusion systems},
\newblock \bibinfo{journal}{Int. J. Numer. Anal. Model.} \bibinfo{volume}{3}
  (\bibinfo{year}{2006}) \bibinfo{pages}{125--140}.
\bibitem[{Higham and Kloeden(2007)}]{HK2007a}
\bibinfo{author}{D.~J. Higham}, \bibinfo{author}{P.~E. Kloeden},
\newblock \bibinfo{title}{Strong convergence rates for backward {E}uler on a
  class of nonlinear jump-diffusion problems},
\newblock \bibinfo{journal}{J. Comput. Appl. Math.} \bibinfo{volume}{205}
  (\bibinfo{year}{2007}) \bibinfo{pages}{949--956}.
\bibitem[{Dareiotis et~al.(2016)Dareiotis, Kumar, and Sabanis}]{DKS2016}
\bibinfo{author}{K.~Dareiotis}, \bibinfo{author}{C.~Kumar},
  \bibinfo{author}{S.~Sabanis},
\newblock \bibinfo{title}{On tamed {E}uler approximations of {SDE}s driven by
  {L}\'evy noise with applications to delay equations},
\newblock \bibinfo{journal}{SIAM J. Numer. Anal.} \bibinfo{volume}{54}
  (\bibinfo{year}{2016}) \bibinfo{pages}{1840--1872}.
\bibitem[{Kumar and Sabanis(2017)}]{KS2017}
\bibinfo{author}{C.~Kumar}, \bibinfo{author}{S.~Sabanis},
\newblock \bibinfo{title}{On tamed {M}ilstein schemes of {SDE}s driven by
  {L}\'evy noise},
\newblock \bibinfo{journal}{Discrete Contin. Dyn. Syst. Ser. B}
  \bibinfo{volume}{22} (\bibinfo{year}{2017}) \bibinfo{pages}{421--463}.
\bibitem[{Hutzenthaler et~al.(2012)Hutzenthaler, Jentzen, and Kloeden}]{Hut02}
\bibinfo{author}{M.~Hutzenthaler}, \bibinfo{author}{A.~Jentzen},
  \bibinfo{author}{P.~E. Kloeden},
\newblock \bibinfo{title}{Strong convergence of an explicit numerical method
  for {SDE}s with nonglobally {L}ipschitz continuous coefficients},
\newblock \bibinfo{journal}{Ann. Appl. Probab.} \bibinfo{volume}{22}
  (\bibinfo{year}{2012}) \bibinfo{pages}{1611--1641}.
\bibitem[{Higham(2011)}]{Hig2011a}
\bibinfo{author}{D.~J. Higham},
\newblock \bibinfo{title}{Stochastic ordinary differential equations in applied
  and computational mathematics},
\newblock \bibinfo{journal}{IMA J. Appl. Math.} \bibinfo{volume}{76}
  (\bibinfo{year}{2011}) \bibinfo{pages}{449--474}.
\bibitem[{Zhang and Ma(2017)}]{ZM2017}
\bibinfo{author}{Z.~Zhang}, \bibinfo{author}{H.~Ma},
\newblock \bibinfo{title}{Order-preserving strong schemes for {SDE}s with
  locally {L}ipschitz coefficients},
\newblock \bibinfo{journal}{Appl. Numer. Math.} \bibinfo{volume}{112}
  (\bibinfo{year}{2017}) \bibinfo{pages}{1--16}.
\bibitem[{Sabanis(2013)}]{Sabanis2013}
\bibinfo{author}{S.~Sabanis},
\newblock \bibinfo{title}{A note on tamed {E}uler approximations},
\newblock \bibinfo{journal}{Electron. Commun. Probab.} \bibinfo{volume}{18}
  (\bibinfo{year}{2013}) \bibinfo{pages}{1--10}.
\bibitem[{Hutzenthaler and Jentzen(2015)}]{HJ2015}
\bibinfo{author}{M.~Hutzenthaler}, \bibinfo{author}{A.~Jentzen},
\newblock \bibinfo{title}{Numerical approximations of stochastic differential
  equations with non-globally {L}ipschitz continuous coefficients},
\newblock \bibinfo{journal}{Mem. Amer. Math. Soc.} \bibinfo{volume}{236}
  (\bibinfo{year}{2015}) \bibinfo{pages}{1--95}.
\bibitem[{Mao(2015)}]{Mao2015}
\bibinfo{author}{X.~Mao},
\newblock \bibinfo{title}{The truncated {E}uler-{M}aruyama method for
  stochastic differential equations},
\newblock \bibinfo{journal}{J. Comput. Appl. Math.} \bibinfo{volume}{290}
  (\bibinfo{year}{2015}) \bibinfo{pages}{370--384}.
\bibitem[{Mao(2016)}]{Mao2016}
\bibinfo{author}{X.~Mao},
\newblock \bibinfo{title}{Convergence rates of the truncated {E}uler-{M}aruyama
  method for stochastic differential equations},
\newblock \bibinfo{journal}{J. Comput. Appl. Math.} \bibinfo{volume}{296}
  (\bibinfo{year}{2016}) \bibinfo{pages}{362--375}.
\bibitem[{{Tan} and {Yuan}(2018)}]{TY2018arXiv}
\bibinfo{author}{L.~{Tan}}, \bibinfo{author}{C.~{Yuan}},
\newblock \bibinfo{title}{{Convergence rates of truncated EM scheme for
  NSDDEs}},
\newblock \bibinfo{journal}{ArXiv e-prints}  (\bibinfo{year}{2018}).
\bibitem[{Kumar and Sabanis(2017)}]{KS2017b}
\bibinfo{author}{C.~Kumar}, \bibinfo{author}{S.~Sabanis},
\newblock \bibinfo{title}{On explicit approximations for {L}\'evy driven {SDE}s
  with super-linear diffusion coefficients},
\newblock \bibinfo{journal}{Electron. J. Probab.} \bibinfo{volume}{22}
  (\bibinfo{year}{2017}) \bibinfo{pages}{1--19}.
\bibitem[{Yang and Wang(2017)}]{YW2017a}
\bibinfo{author}{X.~Yang}, \bibinfo{author}{X.~Wang},
\newblock \bibinfo{title}{A transformed jump-adapted backward {E}uler method
  for jump-extended {CIR} and {CEV} models},
\newblock \bibinfo{journal}{Numer. Algorithms} \bibinfo{volume}{74}
  (\bibinfo{year}{2017}) \bibinfo{pages}{39--57}.
\bibitem[{Przybyl~owicz(2016)}]{Prz2016}
\bibinfo{author}{P.~Przybyl~owicz},
\newblock \bibinfo{title}{Optimal global approximation of stochastic
  differential equations with additive {P}oisson noise},
\newblock \bibinfo{journal}{Numer. Algorithms} \bibinfo{volume}{73}
  (\bibinfo{year}{2016}) \bibinfo{pages}{323--348}.
\bibitem[{Mao et~al.(2016)Mao, You, and Mao}]{MYM2016}
\bibinfo{author}{W.~Mao}, \bibinfo{author}{S.~You}, \bibinfo{author}{X.~Mao},
\newblock \bibinfo{title}{On the asymptotic stability and numerical analysis of
  solutions to nonlinear stochastic differential equations with jumps},
\newblock \bibinfo{journal}{J. Comput. Appl. Math.} \bibinfo{volume}{301}
  (\bibinfo{year}{2016}) \bibinfo{pages}{1--15}.
\bibitem[{Wang and Gan(2010)}]{WG2010}
\bibinfo{author}{X.~Wang}, \bibinfo{author}{S.~Gan},
\newblock \bibinfo{title}{Compensated stochastic theta methods for stochastic
  differential equations with jumps},
\newblock \bibinfo{journal}{Appl. Numer. Math.} \bibinfo{volume}{60}
  (\bibinfo{year}{2010}) \bibinfo{pages}{877--887}.
\bibitem[{Kloeden and Platen(1992)}]{KP1992a}
\bibinfo{author}{P.~E. Kloeden}, \bibinfo{author}{E.~Platen},
  \bibinfo{title}{Numerical Solution of Stochastic Differential Equations},
  \bibinfo{publisher}{Springer-Verlag}, \bibinfo{address}{Berlin},
  \bibinfo{year}{1992}.
\bibitem[{Platen and Bruti-Liberati(2010)}]{Pla2010}
\bibinfo{author}{E.~Platen}, \bibinfo{author}{N.~Bruti-Liberati},
  \bibinfo{title}{Numerical Solution of Stochastic Differential Equations with
  Jumps in Finance}, \bibinfo{publisher}{Springer-Verlag},
  \bibinfo{address}{Berlin}, \bibinfo{year}{2010}.
\bibitem[{Guo et~al.(2017)Guo, Liu, Mao, and Yue}]{Guo2017partial}
\bibinfo{author}{Q.~Guo}, \bibinfo{author}{W.~Liu}, \bibinfo{author}{X.~Mao},
  \bibinfo{author}{R.~Yue},
\newblock \bibinfo{title}{The partially truncated {E}uler-{M}aruyama method and
  its stability and boundedness},
\newblock \bibinfo{journal}{Appl. Numer. Math.} \bibinfo{volume}{115}
  (\bibinfo{year}{2017}) \bibinfo{pages}{235--251}.
\bibitem[{Bao et~al.(2011)Bao, B\"ottcher, Mao, and Yuan}]{Bao2011}
\bibinfo{author}{J.~Bao}, \bibinfo{author}{B.~B\"ottcher},
  \bibinfo{author}{X.~Mao}, \bibinfo{author}{C.~Yuan},
\newblock \bibinfo{title}{Convergence rate of numerical solutions to {SFDE}s
  with jumps},
\newblock \bibinfo{journal}{J. Comput. Appl. Math.} \bibinfo{volume}{236}
  (\bibinfo{year}{2011}) \bibinfo{pages}{119--131}.
\bibitem[{Guo et~al.(2018)Guo, Liu, and Mao}]{Guo2018note}
\bibinfo{author}{Q.~Guo}, \bibinfo{author}{W.~Liu}, \bibinfo{author}{X.~Mao},
\newblock \bibinfo{title}{A note on the partially truncated {E}uler-{M}aruyama
  method},
\newblock \bibinfo{journal}{Applied Numerical Mathematics}
  \bibinfo{volume}{130} (\bibinfo{year}{2018}) \bibinfo{pages}{157--170}.

\end{thebibliography}

\end{document}